\documentclass{amsart}
\usepackage{amssymb,amscd}
\usepackage{graphicx}

\usepackage{latexsym}

\usepackage{amsthm}
\newtheorem{theorem}{Theorem}
\newtheorem*{thm}{Theorem}
\newtheorem{proposition}{Proposition}

\newtheorem{defn}{Definition}
\newtheorem{lemma}{Lemma}
\newtheorem{cor}{Corollary}

\newcommand{\ad}{\rm{ad}}
\newcommand{\bc}{\mathbb{C}}

\newcommand{\bp}{\mathbb{P}}

\newcommand{\br}{\mathbb{R}}

\newcommand{\f}{\mathcal{F}}
\newcommand{\g}{\mathcal{G}}
\newcommand{\q}{\mathcal{Q}}
\newcommand{\y}{\mathcal{Y}}
\newcommand{\z}{\mathcal{Z}}
\newcommand{\tr}{{\rm tr}\hspace{.2mm}}

\newcommand{\ie}{i.e.\ }

\begin{document}

\title{Magnetic monopoles on manifolds with boundary}
\author{Paul Norbury}
\address{\hspace{-.5mm}Department of Mathematics and Statistics\\
University of Melbourne\\Australia 3010.}
\email{pnorbury@ms.unimelb.edu.au}

\date{\today}

\begin{abstract}

\noindent Kapustin and Witten associate a Hecke modification of a holomorphic bundle over a Riemann surface to a singular monopole on a Riemannian surface times an interval satisfying prescribed boundary conditions.  We prove existence and uniqueness of singular monopoles satisfying prescribed boundary conditions for any given Hecke modification data  confirming the underlying geometric invariant theory principle.

\end{abstract}

\maketitle

\section{Introduction}

Dirac monopoles are singular solutions of Maxwell's equation on a Riemannian 3-manifold $Y$.  They arose out of Dirac's study of the quantum theory of electro-magnetism.   Bogomolny-Prasad-Sommerfield, or non-abelian, monopoles are a generalisation of Dirac monopoles to non-abelian theories where the singularities can be smoothed away, \cite{AHiGeo}. For any compact Lie group $G$, they are smooth solutions of the Bogomolny equation 
\begin{equation}  \label{eq:bog}
d^A\Phi=*F_A
\end{equation} 
where $(A,\Phi)$ is a pair given by a $G$-connection $A$ and a section of the adjoint bundle, the Higgs field $\Phi$, defined on a trivial $G$-bundle $E$ over $Y$.  The curvature of $A$ is $F_A$ and $d^A\Phi$ is the covariant derivative of the Higgs field.  The Hodge star $*$ is given with respect to the metric on $Y$.  When $G=U(1)$, the Bogomolny equation (\ref{eq:bog}) reduces to the static Maxwell's equation $d\Phi=*dA$.

On the manifold $Y=\Sigma\times I$ for a Riemannian surface $\Sigma$ and an interval $I=[0,1]$, equipped with the product metric $ds^2=4\rho(z,\bar{z})d\bar{z}dz+dt^2$, the Bogomolny equation (\ref{eq:bog}) uses the full Riemannian metric on $\Sigma$ and splits into a complex equation that depends only on the underlying conformal structure on $\Sigma$ and a real equation:
\begin{equation*}  \label{eq:bogprod}
 [\partial_{\bar{z}}^A,\partial_t^A-i\Phi]=0,\quad  [\partial_{\bar{z}}^A,\partial_z^A]=\rho(z,\bar{z})[\partial_t^A+i\Phi,\partial_t^A-i\Phi].
\tag{1a}
\end{equation*} 
The $A_{\bar{z}}$ part of the connection $A$ determines a holomorphic structure on the restriction $E_t$ of $E$ to each slice $\Sigma\times\{t\}$.  The complex equation in (\ref{eq:bogprod}) can be rewritten
\begin{equation*}  \label{eq:cplxeq}
 \partial_tA_{\bar{z}}=\partial_{\bar{z}}^A(A_t-i\Phi).
\tag{1b}
\end{equation*}
In other words $A_{\bar{z}}$ changes by the complex gauge transformation $\exp[\int_{t_0}^t(A_t-i\Phi)dt]$ so the holomorphic type of the bundle $E_t$ is independent of $t$.  (In fact, (\ref{eq:cplxeq}) holds more generally for a conformal product---when $\rho=\rho(z,\bar{z},t)$ and the {\em conformal} type of the metric on $\Sigma$ is independent of $t$, e.g. the Euclidean metric on $S^2\times\br\cong\br^3-0$.) 

A relaxation of the smoothness of non-abelian monopoles was introduced by Kronheimer over $\br^3$ in \cite{KroMon} and over any closed Riemannian 3-manifold by Pauly in \cite{PauMon}, allowing very special isolated singularities in solutions of (\ref{eq:bog}).  
A singular $G$-monopole $(A,\Phi)$ on $E$ over a 3-manifold $Y$ is a smooth monopole everywhere except at a finite number of points.  In a ball neighbourhood $B^3(y_0)$ of a singularity $y_0\in Y$, there exists a local quotient $\pi:B^4\to B^3(y_0)$ by an action of $U(1)$ that restricts to the Hopf map over 2-spheres foliating $B^3(y_0)$, and sends the single fixed point of the action to $y_0$.  Upstairs, $B^4$ is equipped with a $U(1)$ invariant Riemannian metric that maps to the given metric on $B^3(y_0)\subset Y$, and a bundle $\tilde{E}$ over $B^4$ with a lift of the $U(1)$ action.  The bundle $\tilde{E}$ has a smooth $U(1)$ invariant $G$-connection $\tilde{A}$ satisfying the anti-self-duality (ASD) equation $F_{\tilde A}=-*F_{\tilde A}$.  The Bogomolny equation (\ref{eq:bog}) is a dimensional reduction of the ASD equation, and ${\tilde A}$ induces the monopole $(A,\Phi)$ downstairs which is necessarily singular at $y_0\in Y$.  In particular, given a homomorphism $U(1)\to G$, Dirac monopoles are examples of singular monopoles. 
  
The lift of the $U(1)$ action on the pull-back bundle $\tilde{E}$ upstairs restricts to a homomorphism $\lambda:U(1)\to G$ over the fixed point above $y_0$.   The {\em weight} of the singularity of $(A,\Phi)$ at $y_0$ is defined to be the conjugacy class of $\lambda$.  Conjugacy
classes of homomorphisms $\lambda:U(1)\to G$ are classified by highest weights of the Langlands dual group $^LG$, or equivalently by irreducible representations of $^LG$.

Following Kapustin and Witten \cite{KWiEle} we consider singular monopoles over the product Riemannian manifold $Y=\Sigma\times I$ for a Riemannian surface $\Sigma$ and an interval $I=[0,1]$.  Suppose there is a singular point at $y_0=(z_0,t_0)$.  In this case
(\ref{eq:cplxeq}) shows that $E_-=E_t$ is a well-defined holomorphic bundle for $t<t_0$ then it jumps to give a different holomorphic bundle $E_+=E_t$ for $t>t_0$ which satisfies $E_-\cong E_+$ on $\Sigma-z_0$.  Furthermore, this isomorphism is canonical since if we choose a trivialisation of $E_0$ on $\Sigma-z_0$ then (\ref{eq:cplxeq}) transports this trivialisation along $t\in I$ so it is well-defined on $E_-$ and also on $E_+$.  Thus, a singular monopole on $Y=\Sigma\times I$ gives rise to a Hecke modification (defined below.)  Kapustin and Witten defined the following boundary conditions for singular monopoles on $\Sigma\times I$.
\begin{equation}  \label{eq:bou}
\hspace{-2.6cm} Boundary\ conditions:\quad
A|_{\Sigma\times0}=A_0,\quad\Phi|_{\Sigma\times1}=0
\end{equation}
where $A_0$ is the trivial connection.

Define the gauge group $\g$ to be the space of smooth maps $g:\Sigma\times I\rightarrow G$ satisfying $g|_{\Sigma\times0}=I$.  It acts on pairs by $g\cdot(A,\Phi)=(gAg^{-1}-dg\cdot g^{-1},g\Phi g^{-1})$ and two solutions are gauge equivalent if they are related by a gauge transformation.  
\begin{defn}
The moduli space of monopoles $\z(\lambda,y_0)$ is the space of gauge equivalence classes of solutions of (\ref{eq:bog}) over $Y=\Sigma\times I$, with singularity at $y_0\in Y$ of weight $\lambda$ satisfying the boundary conditions (\ref{eq:bou}).
\end{defn} 

Note that there is an action of $G$ on $\z(\lambda,y_0)$ given by constant gauge transformations.  Such gauge transformations preserve the trivial connection $A_0$ and give gauge inequivalent monopoles since they violate the condition $g|_{\Sigma\times0}=I$. 

The surface $\Sigma$ possesses a Riemannian metric.  In the following two definitions, we consider only its conformal class.  Let $G^c$ be the complexification of $G$ and denote a holomorphic $G$-bundle to mean a holomorphic $G^c$-bundle. 

\begin{defn}
A Hecke modification of a holomorphic $G$-bundle $E_-$ over a Riemann surface $\Sigma$ at $z_0\in\Sigma$, is a pair $(E_+,h)$ consisting of a holomorphic $G$-bundle $E_+$ and an isomorphism
\[ h:E_-|_{\Sigma-z_0}\stackrel{\cong}{\longrightarrow} E_+|_{\Sigma-z_0}.\]
\end{defn}
Equivalently, a Hecke modification of a holomorphic bundle $E_-$ at $z_0\in\Sigma$ is an exact sequence of sheaves
\[0\to E_-\to E_+\to\q\to 0\]
where $E_-$ and $E_+$ are locally free and $\q$ is a skyscraper sheaf at $z_0$ (and we may have to replace $E_+$ with $E_+\otimes\mathcal{O}(Nz_0)$ so that the isomorphism outside $z_0$ extends.)

With respect to a local chart on $\Sigma$ the bundles $E_-$ and $E_+$ can be defined using the same transition functions (representatives of a Cech cohomology class) except in a neighbourhood $U$ of $z_0$.  Choose a local parameter $z$ in $U$ and let $\lambda:U(1)\to G$ be a group homomorphism.  The Hecke modification $(E_+,h)$ of $E_-$ is of weight $\lambda$ if local trivialisations can be chosen so that for any local chart $V$ overlapping $U$ the transition functions for $E_+$ and $E_-$ are related by
\[ g^+_{UV}(z)=\lambda(z)g^-_{UV}(z).\]
\begin{defn}
Define $\y(E_-,\lambda,z_0)$ to be the space of Hecke modifications of $E_-$, \ie the space of all pairs $(E_+,h)$ with $h:E_-\stackrel{\cong}{\to} E_+$ on $\Sigma-\{z_0\}$, of weight $\lambda$.  Suppress $E_-$ when it is the trivial bundle and simply write $\y(\lambda,z_0)$.  
\end{defn}
The space $\y(\lambda,z_0)$ depends only on the conjugacy class of $\lambda$ which is called the {\em weight} of the Hecke modification.  See \cite{FBeVer} for further details on Hecke modifications and applications.

Kapustin and Witten \cite{KWiEle} showed that with the boundary conditions (\ref{eq:bou}) and the condition that gauge transformations are trivial on $\Sigma\times0$, the second equation in (\ref{eq:bogprod}) is the zero set of the moment map of the action of the gauge group $\g$ on the space of pairs $(A,\Phi)$ on the bundle $E$ with respect to a natural symplectic structure.   The Hecke modification associated to a singular monopole on $\Sigma\times I$ uses only the first equation in (\ref{eq:bogprod}) which is preserved by the group $\g^c$ of $G^c$ gauge transformations.   The symplectic quotient of pairs $(A,\Phi)$ satisfying the first equation in (\ref{eq:bogprod}) is the moduli space of monopoles $\z(\lambda,y_0)$.  The forgetful map from solutions of (\ref{eq:bogprod}) to solutions of the first equation in (\ref{eq:bogprod}) defines
\[ \f:\z(\lambda,y_0)\to\y(\lambda,z_0).\] 
(The fact that the weight $\lambda$ is respected requires proof.  See Proposition~\ref{th:imf}.)

Given a $\g$ action on a symplectic manifold $(M,\omega)$ with moment map $\mu$, if the action of $\g$ lifts to a complexified action of $\g^c$, the basic principle of geometric invariant theory is a correspondence
\[ \mu^{-1}(0)/\g\cong M^s/\g^c\]
for a subset $M^s\subset M$ of stable points.
Such a correspondence is an {\em existence} and {\em uniqueness} result: in each stable $\g^c$ orbit there exists a unique $\g$ orbit that maps to $0$ under $\mu$.  This principle applied to the symplectic manifold of pairs $(A,\Phi)$ satisfying the first equation in (\ref{eq:bogprod}) which is invariant under the action of the complexified gauge group $\g^c$ predicts that $\f:\z(\lambda,y_0)\to\y(\lambda,z_0)$ is a homeomorphism between the symplectic quotient $\z(\lambda,y_0)$ and the quotient by complex gauge transformations $\y(\lambda,z_0)$.  The main result of this paper is this existence and uniqueness result for solutions of the boundary value Bogomolny equation.  We show that to each Hecke modification there corresponds a monopole and that such a monopole is unique, yielding the following theorem.

\begin{theorem}  \label{th:hom}
For any conjugacy class of homomorphisms from the circle to a compact Lie group $\lambda:U(1)\to G$, the map $\f:\z(\lambda,y_0)\to\y(\lambda,z_0)$ is a homeomorphism for $y_0=(z_0,t_0)\in\Sigma\times I$.
\end{theorem}

More generally, for a collection of distinct points $y_1,...,y_k\in\Sigma\times I$ denote the moduli space of  monopoles with prescribed singularities of weights $\lambda_1,...,\lambda_k$ at the $y_i$ by $\z(\lambda_1,y_1,...,\lambda_k,y_k)$.   Theorem~\ref{th:hom} is a special case of:
\begin{theorem}  \label{th:main}
For any collection of conjugacy classes of homomorphisms from the circle to a compact Lie group $\lambda_i:U(1)\to G$, when the $z_i$ are distinct the map 
\[ \f:\z(\lambda_1,y_1,...,\lambda_k,y_k)\to\y(\lambda_1,z_1)\times...\times\y(\lambda_k,z_k)\]
is a homeomorphism for $y_i=(z_i,t_i)\in\Sigma\times I$.
\end{theorem}
A version of the theorem still holds when the $z_i$ are not distinct.  Suppose, for example, that $z_1=z_ 2$ and $t_1<t_2$.  Then Hecke modifications of weight $\lambda_1$ of the trivial bundle still give rise to $\y(\lambda_1,z_1)$ whereas the Hecke modifications of weight $\lambda_2$ act not on the trivial bundle, instead on the pairs $(E_+,h)\in\y(\lambda_1,z_1)$.  For each $(E_+,h)$ the space of Hecke modifications of weight $\lambda_2$ is isomorphic to $\y(\lambda_2,z_1)$ non-canonically.  In other words, one ends up with a $\y(\lambda_2,z_1)$ fibre bundle over $\y(\lambda_1,z_1)$.  Thus we can replace the product in Theorem~\ref{th:main} by appropriate fibre bundles and the theorem goes through.

For concreteness we give here a particular case of the theorem.  Let $G=PSU(2)$, and represent each weight $\lambda=(m,-m)$ by the half-integer $m$.  Consider any moduli space of $PSU(2)$ monopoles where all singularities have weight $1/2$.  Such a moduli space is compact because bubbling---resulting from a non-compact sequence of monopoles---can occur only at singularities, however the weight $1/2$ singularity is too small to allow this, see \cite{KWiEle}.  There is a two-sphere of Dirac monopoles inside the moduli space, corresponding to homomorphisms $U(1)\to PSU(2)$ of weight $(1/2,-1/2)$.  It is easy to show that $\y(1/2,z_1)\cong\bp^1$ so a consequence of Theorem~\ref{th:hom} is that {\em all} monopoles in $\z(1/2,p)$ are Dirac monopoles.  In contrast, for
\[\z(1/2,y_1,...,1/2,y_k)\cong\y(1/2,z_1)\times...\times\y(1/2,z_k)\]  
the diagonal $\bp^1\stackrel{\Delta}{\hookrightarrow}(\bp^1)^k$ consists of Dirac monopoles and the other monopoles are truly nonabelian.  The idea of the proof of the existence of these nonabelian monopoles is to show that the continuous and differentiable function $\f$ is an open mapping, which is a linear issue.  Using this and the compactness of the domain, the method of continuity applies to give a monopole for each $k$-tuple of Hecke modifications in $\y(1/2,z_1)\times...\times\y(1/2,z_k)$.  More generally, when the domain is not compact, an understanding of how compactness fails, \ie bubbling, allows one to prove $\f$ has closed image.  

By lifting the problem to four dimensions, we give a second indirect proof of the existence and uniqueness of singular monopoles using a result of Donaldson and Simpson \cite{DonBou,SimCon} for anti-self-dual connections on a holomorphic bundle over a K\"ahler surface with boundary.  This viewpoint is interesting even before applying the gauge theory results of  \cite{DonBou,SimCon} since it produces a single holomorphic bundle over a complex surface to encode a Hecke modification which consists of {\em two} holomorphic bundles over a curve.

Theorems~\ref{th:hom} and \ref{th:main} are analogous to a theorem of Donaldson \cite{DonNah} which identifies the moduli space of charge $k$ $SU(2)$ monopoles over $\br^3$ with the space of degree $k$ rational maps from $S^2$ to $S^2$.  The methods in this paper are different to the Nahm equation methods used by Donaldson in  \cite{DonNah} and more closely related to the methods introduced by Donaldson for the study of instantons over algebraic surfaces \cite{DonAnt,DonBou} and applied to monopoles in \cite{JarCon,JNoCom}.

Inside the space $\y(E_-,\lambda,z_0)$ is the subvariety of pairs $(E_+,h)$ with the holomorphic bundle $E_+$ fixed.  Kapustin and Witten \cite{KWiEle} showed that when $E_+$ is a stable holomorphic bundle this corresponds to singular monopoles with new boundary conditions: the restriction of the connection $A$ to each end $\Sigma\times 0$ and $\Sigma\times 1$ is a flat connection.  The techniques in this paper can be used to give a proof that this gives a 1-1 correspondence.

{\bf Acknowledgements:}  The author is grateful to Peter Kronheimer, Mark Stern and Benoit Charbonneau for useful conversations.

\section{The map $\f$}   \label{sec:hecke}
The map $\f:\z(\lambda,y_0)\to\y(\lambda,z_0)$ is defined in two steps.  Firstly it is defined as a map $\f:\z(\lambda,y_0)\to\Omega G$, where $\Omega G$ is the loop group of based smooth maps from $S^1$ to $G$.  The space of Hecke modifications of weight $\lambda$ can be identified with a subvariety of the loop group $\y(\lambda,z_0)\subset\Omega G$.  Thus, secondly one needs to show that the image of $\f$ is contained in  $\y(\lambda,z_0)$.  Kapustin and Witten \cite{KWiEle} used a topological argument to prove the second step when $G=U(1)$, and using this proved it for $G$ monopoles which are sums of Dirac monopoles.  

In this section, given a monopole $(A,\Phi)$, we give an explicit construction of $\gamma\in\Omega G$, proving continuity properties. Furthermore we show that the image of $\f$ is indeed contained in $\y(\lambda,z_0)\subset\Omega G$.  The construction of $\gamma=\f(A,\Phi)$ is purely local so it generalises immediately to many singularities $(\gamma_1,...,\gamma_k)=\f(A,\Phi)$ for $(A,\Phi)\in\z(\lambda_1,y_1,...,\lambda_k,y_k)$. We begin with definitions of the loop group and the space of Hecke modifications.

\subsection{The loop group}
Let $LG^c$ be the loop group of smooth maps from the circle $\{|z|=1\}\subset\bc$ to $G^c$ and $L^+G^c\subset LG^c$ those smooth maps that extend from the circle to holomorphic maps of the disk to $G^c$.  The quotient $LG^c/L^+G^c$ of $LG^c$ by the right action of $L^+G^c$ behaves like an infinite dimensional Grassmannian.  In particular there is an isomorphism $LG^c/L^+G^c\cong\Omega G$ to the loop group of based smooth maps from the circle to $G$.   Similarly define $L{\bf g}^c/L^+{\bf g}^c$ for the Lie algebra ${\bf g}^c$.  If $\lambda:S^1\to G^c$ is a group homomorphism, then the orbit $L^+G^c\cdot\lambda\subset LG^c/L^+G^c$ induced by the left action of $L^+G^c$ on $LG^c$ consists of loops (with coset representatives) that are algebraic, or equivalently extend to maps $\bc^*\to G^c$, and the union of such orbits is denoted by $\Omega_{\rm pol}G\subset LG^c/L^+G^c$.    All of these facts can be found in \cite{PSeLoo}.

\begin{lemma}
The group $\Omega_{\rm pol}G$ and its strata $L^+G^c\cdot\lambda$, indexed by conjugacy classes of homomorphisms $\lambda:S^1\to G^c$, are independent of the coordinate $z$ so they can be defined using arbitrary local coordinates on a Riemann surface. 
\end{lemma}
\begin{proof}
Given a local coordinate $z$, choose a new local coordinate $z(w)$ in the disk $|z|<1$ such that $z(0)=0$ and for the moment assume that $|w|<1$ implies $|z|<1$.  (More generally, any two disks $|z|<1$ and $|w|<1$ contain a smaller disk in their intersection so we apply this reasoning twice.)  The restriction map $\gamma(z)\mapsto\gamma(|w|=1)$ is a group homomorphism $F:LG^c\to LG^c$, \ie $F(\gamma_1\gamma_2)=F(\gamma_1)F(\gamma_2)$, which induces a one-to-one map $\Omega_{\rm pol}G\to\Omega_{\rm pol}G$, and preserves the stratification by $L^+G^c\cdot\lambda$ for any homomorphism $\lambda$, due to the following properties:
\begin{enumerate}
\item[(i)] $F$ maps $L^+G^c$ to $L^+G^c$.
\item[(ii)] $F$ maps $\lambda\cdot L^+G^c$ to $\lambda\cdot L^+G^c$.
\item[(iii)] The induced map $L^+G^c\cdot\lambda\subset\Omega_{\rm pol}G$ to $L^+G^c\cdot\lambda\subset\Omega_{\rm pol}G$ is surjective
\end{enumerate}
Property (i) is immediate.
To prove (ii), put $\lambda(z)={\rm diag}(z^{m_1},...,z^{m_n})$, then 
\[ F(\lambda)={\rm diag}(z(w)^{m_1},...,z(w)^{m_n})\equiv{\rm diag}(w^{m_1},...,w^{m_n})({\rm mod}\ L^+G^c)\] 
since $z(w)/w:\{|w|\leq 1\}\to G^c$ is holomorphic.  Surjectivity in (iii) uses the fact that $\gamma=p(w)\lambda(w)\in L^+G^c\cdot\lambda({\rm mod}\ L^+G^c)$ extends to $\gamma:\bc^*\to G^c$ so we may assume $p(w)$ extends, and by restriction $\gamma(w(z))=p(w(z))\lambda(w(z))$.
\end{proof}
Given a point $z_0$ on a Riemann surface $\Sigma$ and $\gamma:D\backslash\{ z_0\}\to G^c$ a holomorphic map  from a deleted neighbourhood of $z_0$ to $G^c$ choose any local coordinate $z$, denote by $D$ the disk $|z|<\epsilon$ and rescale $z$ so that $\epsilon=1$.  This gives an element of the loop group $\gamma:\{|z|=1\}\to G^c$ and any other choice of local coordinate gives another element in the loop group which is the image of $\gamma$ under a group homomorphism.

\subsection{The space of Hecke modifications}
Elements of the loop group $\Omega G$ para\-met\-rise pairs $(E,h)$ consisting of a holomorphic bundle $E$ over a curve trivialised by $h$ over the complement of a point.  Thus the space of Hecke modifications of a holomorphic bundle over a curve at a point is identified as a subvariety of the loop group $\Omega G$ of smooth loops, known as the affine grassmannian \cite{FBeVer} of algebraic loops.  This is a convenient description which produces, for example,  a calculation of the dimension of $\y(E,z_0,\lambda)$.

\begin{lemma}   \label{th:yinlg}
For the trivial holomorphic $G^c$ bundle $E$ over $\Sigma$,
\[
\y(E,z_0,\lambda)\cong L^+G^c\cdot\lambda\subset LG^c/ L^+G^c.
\]
\end{lemma}
\begin{proof}  
For any holomorphic $G^c$ bundle over $\Sigma$ defined via  transition functions over an open cover of $\Sigma$, if $G$ is connected it is sufficient to use the two charts
\[\Sigma=\left(\Sigma-z_0\right)\cup D_{z_0}\] 
for $D_{z_0}$ a disk neighbourhood of $z_0$ in $\Sigma$.  This is because $\Sigma-z_0$ is a Stein manifold so topological and holomorphic classification of bundles coincide \cite{GraHol}.  In particular, $\Sigma-z_0$ retracts to a one-dimensional spine over which any principal $G$ bundle with connected $G$ is trivial.  Hence, a holomorphic $G^c$ bundle is trivial over $\Sigma-z_0$ and $D_{z_0}$.  Its transition function is a map $\gamma:D_{z_0}\cap(\Sigma-z_0)\to G^c$ which we identify with $\gamma:S^1\to G^c$ or equivalently $\gamma\in LG^c$.  The bundle is also defined by the equivalent transition function $p_-\gamma p_+$ for $p_{\pm}\in LG^c$ that extend to holomorphic maps $p_-:\Sigma-z_0\to G^c$ and $p_+:D_{z_0}\to G^c$, and we write $p_+\in L^+G^c$.

A bundle $E_+$ is defined by $\gamma:S^1\to G^c$ well-defined up to $p_-\gamma p_+$ so a pair $(E_+,h)\in\y(E_-,z_0,\lambda)$ is defined by $\gamma\in LG^c$ well-defined up to $\gamma\mapsto\gamma p_+$ for $p_+\in L^+G^c$, or equivalently $\gamma\in LG^c/L^+G^c\cong\Omega G$.  The pair is a bundle equipped with an isomorphism $h:E_+|_{\Sigma-z_0}\to E_-|_{\Sigma-z_0}$ and $h$ serves to fix the freedom in the transition function $\gamma$ on the left.  With respect to the same trivialisation, $E_-$ is defined by $q$, say.  As defined in the introduction, for a group homomorphism $\lambda:S^1\to G^c$ a Hecke modification of $E_-$ of weight $\lambda$ replaces $q$ by $\gamma=q\lambda$.  We assume that $E_-$ is the trivial bundle, so $q\in L^+G^c$ and $\gamma\in L^+G^c\cdot\lambda$.  Moreover, $\y(\lambda)\cong L^+G^c\cdot\lambda$.
\end{proof}

An immediate application of Lemma~\ref{th:yinlg} is a calculation of the dimension of $\y(E,z_0,\lambda)$.  For any homomorphism $\lambda:S^1\to G^c$, the dimension of its orbit $L^+G^c\cdot\lambda\subset LG^c/ L^+G^c$ follows from understanding the isotropy subgroup 
\[ I_{\lambda}=\{ \gamma\in L^+G^c\ |\ \gamma\cdot\lambda\equiv\lambda({\rm mod\ }L^+G^c)\}\]
which consists of $\gamma\in L^+G^c$ such that $\lambda^{-1}\gamma\lambda\in L^+G^c$.  This is most easily understood on a Lie algebra level.  For example, when $G=SU(2)$ or $PSU(2)$ and $\lambda={\rm diag}(z^m,z^{-m})$ then the Lie algebra ${\bf i}_{\lambda}$ of $I_{\lambda}$ consists of $\xi\in L^+{\bf g}^c$ satisfying $\lambda^{-1}\xi\lambda\in L^+{\bf g}^c$.  Equivalently, the upper diagonal element $\xi_{12}$ of $\xi$ is divisible by $z^{2m}$.  The tangent space of $L^+G^c\cdot\lambda$ is identified with the quotient $L^+{\bf g}^c/{\bf i}_{\lambda}$ which consists of $\eta\in L^+{\bf g}^c$ such that $\eta_{jk}=0$ for $jk\neq(12)$ and $\eta_{12}=p(z)$ with degree $p(z)<2m$.  Thus 
\begin{equation}   \label{eq:dimy}
\dim\y(E,z,m)=2m.
\end{equation} 
The argument easily generalises to higher rank groups.

\subsection{Definition of $\f$}   \label{sec:defnf}
Given a singular monopole $(A,\Phi)$ defined on a trivial $G$ bundle over $Y=\Sigma\times I$ with singularities at $y_i=(z_i,t_i)\in\Sigma\times(0,1)$, $i=1,...,k$ we explicitly define $\f(A,\Phi)$ using the ODE 
\begin{equation}   \label{eq:scat}
(\partial_t^A-i\Phi)g=0
\end{equation} 
with $g(z,0)=I$ defined along $\{z\}\times[0,1]$ for each $z\in\Sigma$.  Existence of solutions of ODEs guarantees $g(z,t)$ exists.  Strictly, to make sense of (\ref{eq:scat}) one needs to multiply it by $g^{-1}$ on the left so that the expression lives in the Lie algebra, and (\ref{eq:scat}) says that $g$ is a complex gauge transformation that sends $A_t-i\Phi$ to 0.  

Near each singularity, the solution $g(z,t)$ of (\ref{eq:scat}) is defined only on $z_i\times[0,t_i)$ and does not extend to $z_i\times(t_i,1]$.  One can find a solution of (\ref{eq:scat}) continuous in a neighbourhood of $z_i\times(t_i,1]$ necessarily given by $g(z,t)\gamma_i(z)$ where $\gamma_i(z)$ is independent of $t$ and defined in a deleted neighbourhood of $z_i$.  
\begin{proposition}   \label{th:definef}
Given a solution $g(z,t)$ of (\ref{eq:scat}) satisfying $g(z,0)=I$ there exists a holomorphic map 
\[\gamma_i:U_i-\{z_i\}\to G^c\] 
where $U_i\subset\Sigma$ is a neighbourhood of $z_i$, unique up to $\gamma_i(z)\mapsto\gamma_i(z)p(z)$ for holomorphic $p:U\to G^c$, satisfying $g(z,t)\gamma_i(z):U_i\times(t_i,1]\to G^c$ solves (\ref{eq:scat}) for each $z\in U_i$.  The map $\f_i(A,\Phi)=\gamma_i$ defines a continuous and differentiable map 
\[\f_i:\z(\lambda_1,y_1,...,\lambda_k,y_k)\to\Omega_{\rm pol}G\subset LG^c/L^+G^c.\]
\end{proposition}
\begin{proof}
The solutions $g(z,t)$ are continuous and differentiable in $z$ since they solve ODEs on finite intervals with continuous and differentiable coefficients.   

Choose a local coordinate on $\Sigma$ with $z=0$  corresponding to $z_i$.  At $z=0$ the solution $g(0,t)$ does not extend past $t=t_i$ \ie $g$ is defined on $\Sigma\times I-\{z_i\}\times(t_i,1]$.  For any $h(z)$ holomorphic in $z\in U\subset\Sigma$, $g(z,t)h(z)$ is a local solution of (\ref{eq:scat}).  For $U$ a disk neighbourhood of $z_i\in\Sigma$ choose $\gamma_i(z)$ so that $g(z,t)\gamma_i(z)$ is well-defined on $U\times I-\{z_i\}\times[0,t_i)$.  

The existence of $\gamma_i(z)$ is a consequence of the existence of solutions to (\ref{eq:scat}) that are also covariantly holomorphic as follows.   The solution $g(z,t)$ of (\ref{eq:scat}) satisfying $g(z,0)=I$ also satisfies $\partial_{\bar{z}}^Ag=0$.  This is a consequence of $[\partial_{\bar{z}}^A,\partial_t^A-i\Phi]=0$, the first equation of (\ref{eq:bogprod}), which implies that $\partial_{\bar{z}}^Ag$ is also a solution of (\ref{eq:scat}).  Thus by uniqueness of solutions, $\partial_{\bar{z}}^Ag=g\mu(z)$ where $\mu(z)$ is independent of $t$.  Since $A$ is trivial along $\Sigma\times0$, $\partial_{\bar{z}}^Ag|_{\Sigma\times 0}=\partial_{\bar{z}}g|_{\Sigma\times 0}=0$ hence $\mu(z)=0$ and $\partial_{\bar{z}}^Ag=0$. Now, let $g_i(z,t)$ be any solution of (\ref{eq:scat}) on $U\times I-\{z_i\}\times[0,t_i)$, satisfying $\partial_{\bar{z}}^Ag_i(z,t)=0$.  By uniqueness of solutions, $g_i(z,t)=g(z,t)\gamma_i(z)$ for some $\gamma_i(z)$ holomorphic in $z\in\bc^*$ since $\partial_{\bar{z}}^Ag(z,t)=0=\partial_{\bar{z}}^Ag_i(z,t)$ proving that $\gamma_i(z)$ exists.

Note that $\gamma_i(z)$ can be replaced by $\gamma_i(z)p(z)$ for any holomorphic map $p:U\to G^c$.  In terms of the loop group, by restricting to the circle $\partial U$, $\gamma_i(z)\in LG^c$ is equivalent to $\gamma_i(z)p(z)$ $({\rm mod}\ L^+G^c)$.  Thus $\f_i(A,\Phi)=\gamma_i(z)\in LG^c/L^+G^c$.  It is continuous and differentiable since $g$ depends continuously and differentiably on $(A,\Phi)$.

The space of smooth loops $LG^c/L^+G^c$ is larger than $\Omega_{\rm pol}G$, the union of orbits $L^+G^c\cdot\lambda$ over all homomorphisms $\lambda$, \cite{PSeLoo}.  {\em A priori} $\gamma_i=\f_i(A,\Phi)$ is a smooth loop not necessarily contained in an orbit $L^+G^c\cdot\lambda$.  Using the Kato inequality
\[|\partial_t|g||\leq|\partial_t^Ag|=|\Phi g|\quad\leq\quad|\Phi||g|\]
hence
\[|\partial_t\ln|g||\leq|\Phi|\quad\leq\quad\frac{c}{r}\]
where $r$ is the distance from $y_i$.  In particular, $g$ has at worst algebraic singularities so $\gamma_i$ has at worst poles and $\gamma_i=\f_i(A,\Phi)\in\Omega_{\rm pol}G$.
\end{proof}

From Proposition~\ref{th:definef}, $(A,\Phi)\in\z(\lambda_1,y_1,...,\lambda_k,y_k)$ maps to $\f_i(A,\Phi)=\gamma_i\in\y(\lambda_i',z_0)$ for some $\lambda'_i$, or equivalently $\gamma_i\in L^+G^c\cdot\lambda_i'\subset LG^c/L^+G^c$.   We need to show that $\lambda_i=\lambda'_i$.  In order to do this, we find a good gauge for the equation (\ref{eq:scat}).   

\subsection{A good gauge.}   \label{sec:gauge}
As described in the introduction, a singular monopole in a ball neighbourhood $B^3(y_i)$ of a singularity $y_i\in Y$ comes from the dimensional reduction of a local quotient $\pi:B^4\to B^3(y_i)$ by an action of $U(1)$.  Upstairs, $B^4$ is equipped with a $U(1)$ invariant Riemannian metric $g^4$ conformally equivalent to the pull-back of the given metric $g^3$ on $B^3(y_i)\subset Y$ together with a metric in the circle direction, which we will see is determined by a positive harmonic function $V$ on $Y-\{y_1,y_2,...,y_k\}$.  The pull-back bundle upstairs $\tilde{E}$ is equipped with a lift of the $U(1)$ action and a smooth $U(1)$ invariant connection $\tilde{A}$ satisfying the ASD equation $F_{\tilde A}=-*F_{\tilde A}$.  We need to be more explicit about the metric upstairs.  This is described in \cite{KroMon,PauMon}.  A metric $g^4$ on a $U(1)$ bundle over a 3-manifold, in this case $\pi:B^4\to B^3(y_i)$, arises from a $U(1)$ singular monopole, which can be encoded by a positive harmonic function $V$ on $B^3(y_i)$ that is singular at $y_i$ and satisfies $rV\to1/2$ as $r\to 0$.  Let $d\theta+\omega$ be the $U(1)$ 
singular monopole, equivalently a $U(1)$ invariant 1-form on $B^4-0$, satisfying $*d\omega=dV$.  Then the monopole $(A,\Phi)$ comes from \[\tilde{A}=\pi^*A-V^{-1}\Phi(d\theta+\omega)\] which is ASD with respect to the metric on $B^4$ given by 
\[g^4=V\pi^*g^3+V^{-1}(d\theta+\omega)^2\] 
(where we have abused terminology by identifying $V$ with $\pi^*V$ and $\Phi$ with $\pi^*\Phi$.)

Near a singular point $y_i$, upstairs in the local model over $B^4$ choose a $U(1)$-equivariant gauge.  The $\theta$ component of the invariant connection $\tilde{A}(\partial/\partial\theta)$ must vanish at the fixed point in an equivariant gauge.  Thus, the weight $\lambda$ appears in a $U(1)$-invariant gauge since it is related to an equivariant gauge by $g$ satisfying $g^{-1}\partial_{\theta}g\to\lambda$ as $r\to 0$, and hence in this gauge $V^{-1}\Phi\to\lambda$ as $r\to 0$.  Furthermore, the ASD equation gives $|F_{\tilde{A}}|^2=2|d^A(V^{-1}\Phi)|^2$ and since $|F_{\tilde{A}}|$ is bounded, $|d^A(V^{-1}\Phi)|$ is bounded.  To summarise, and using the fact that $V\sim 1/2r$, for $r$ the distance from $y_i$, a singular monopole $(A,\Phi)$ with singularity of weight $\lambda\in{\bf g}$ at $y_i$ satisfies
\begin{equation}  \label{eq:sing} 
r\Phi\rightarrow\frac{1}{2}\lambda,\quad |d^A(r\Phi)|\ {\rm is\ bounded\ as\ } r\rightarrow 0.
\end{equation}

We can find a gauge for $A$ so that the growth (\ref{eq:sing}) controls the growth of $A_t-i\Phi$ which controls growth of solutions of (\ref{eq:scat}) as follows.  Near a singularity, choose a gauge in which $\Phi$ is diagonal.  There is a maximal torus $T\subset G$ (with Lie algebra ${\bf t}$) of gauge transformations that preserve $\Phi$ being diagonal, and we can use these to gauge away the ${\bf t}$ component of $A_t$.  Equivalently, $\Phi\in{\bf t}$ and $A_t\in{\bf t}^{\perp}$.  Furthermore, for $\lambda={\rm diag}(m_1,...,m_n)$ if some of the $m_j$ coincide, this gives rise to a larger gauge group preserving $\lambda$, that can be used to gauge away the component of $A_t$ in this direction.  In this gauge, 
\[
\Phi=\frac{1}{r}\left(\begin{array}{ccc}i\phi_1&&0\\&\ddots&\\0&&i\phi_n\end{array}\right),\quad
\partial^A_t=\frac{\partial}{\partial t}+\left(\begin{array}{ccc}0&&\alpha\\&\ddots&\\-\alpha^*&&0\end{array}\right)
\]
where $\alpha_{jk}\to0$ if $m_j=m_k$.  The $\phi_j$ satisfy $\lim_{r\to 0}\phi_j=m_j$.  Furthermore, $|\alpha|$ is bounded since as described above $|d^A(V^{-1}\Phi)|$ is bounded, so $|\partial^A_t(V^{-1}\Phi)|$ is bounded and hence so are its off-diagonal components $(\phi_j-\phi_i)\alpha_{jk}$.  When $m_k-m_j\neq 0$ this implies that $\alpha_{jk}$ is bounded, and when $m_k-m_j=0$, $\alpha_{jk}\to 0$ and in particular it is bounded.

\begin{proposition}  \label{th:imf}
The image of $\f_i$ lies in $\y(\lambda_i,z_i)$.
\end{proposition}
\begin{proof}
Proposition~\ref{th:definef} shows that a solution $g(z,t)$ of $(\partial^A_t-i\Phi)g=0$ with $g(0)=I$ uniquely defines $\gamma_i(z)\in LG^c/L^+G^c$ via the property that $g(z,t)\gamma_i(z)$ extends to $\{z_i\}\times(t_i,1]$.  Furthermore, it is proven that $\gamma_i(z)$ is an algebraic loop, or in other words there exists a homomorphism $\lambda_i'(z)$ such that $\gamma_i(z)=p_i(z)\lambda_i'(z)\in LG^c/L^+G^c$ where $p_i(z)\in L^+G^c$ maps a small disk holomorphically to $G^c$.  By replacing $g(z,t)$ by $g(z,t)p_i(z)$ we may assume $\gamma_i(z)=\lambda_i'(z)$ so in particular it is diagonal.  In the gauge described above (\ref{eq:scat}) becomes 
\[\partial_tg=-(A_t-i\Phi)g=\left(\begin{array}{ccc}-\phi_1/r&&-\alpha\\&\ddots&\\\alpha^*&&-\phi_n/r\end{array}\right)g\]
where $|\alpha|$ is bounded.  Choose a small ball of radius $\delta>0$ around the singularity.  By (\ref{eq:sing}) we can choose $\delta$ small enough so that
\[ m_j-\epsilon<2\phi_j<m_j+\epsilon,\quad r<\delta.\]
Integrating the ODE yields $g(\delta)=\exp(M)g(-\delta)$ for
\[ M=\int^{\delta}_{-\delta}\left(\begin{array}{ccc}-\phi_1/r&&-\alpha\\&\ddots&\\\alpha^*&&-\phi_n/r\end{array}\right)dt=\Lambda+M_0\]
where $M_0$ is bounded and $\Lambda$ is diagonal with components $-\int^{\delta}_{-\delta}\phi_j/rdt$.
\begin{lemma}
\[\int^{\delta}_{-\delta}\frac{c}{r}dt=-c\ln|z|^2+B,\quad B{\rm\ bounded.}\]
\end{lemma}
\begin{proof}
The metric $ds^2=4\rho(z,\bar{z})d\bar{z}dz+dt^2$ implies that in a small neighbourhood of $(z,t)=(0,t_i)$ there exists positive constants $c_1$ and $c_2$ such that
\[ c_1|z|^2+t^2=r_1^2<r^2<r_2^2=c_2|z|^2+t^2\]
where we have changed coordinates $t\mapsto t-t_i$.  Thus 
\[\int^{\delta}_{-\delta}\frac{1}{r}dt<\int^{\delta}_{-\delta}\frac{1}{r_1}dt=\ln\frac{\delta+r_1(\delta)}{-\delta+r_1(\delta)}=\ln
\frac{(\delta+r_1(\delta))^2}{c_1|z|^2}
=-\ln|z|^2+B_1\]
where $B_1$ is bounded and $r_1(\delta)$ is $r_1$ evaluated at $t=\delta$.  We can reverse the inequality by replacing $r_1$ with $r_2$, and the lemma follows.
\end{proof}
By increasing $\epsilon$ slightly we can remove the bounded term and
\[ (m_j+\epsilon)\ln|z|<-\int^{\delta}_{-\delta}\frac{\phi_j}{r}dt<(m_j-\epsilon)\ln|z|\]
so $\Lambda$ has components lying between $(m_j\pm\epsilon)\ln|z|$.

The matrix $M$ is diagonalisable and its eigenspaces tend uniformly to the standard basis as $r\to 0$ with eigenvalues lying between $(m_j\pm\epsilon)\ln|z|$.  Thus $\exp(M)$ has eigenvalues lying between $|z|^{m_j\pm\epsilon}$ as $r\to 0$ and eigenspaces tending uniformly to the standard basis.  Since $g(z)\lambda_i'(z)$ is bounded, $\exp(M)\lambda_i'(z)$ is bounded, and  thus $\lambda_i'(z)=\lambda_i(z)$.
\end{proof}

The proof of Proposition~\ref{th:imf} also gives the behaviour of the growth of $g(z,t)$ near each singularity $y_i$.  We will abuse notation and use $r$ to denote the distance from each $y_i$.
\begin{proposition}  \label{th:imfc}
If $\f_i(A,\Phi)=\gamma_i(z)=p_i(z)\lambda_i(z)$ then near $y_i$ as $r\to 0$ the solution of (\ref{eq:scat}) is asymptotic to
\[g(z,t)\sim\exp\Lambda(z,t)\cdot p_i(z)^{-1}\] 
where $\exp\Lambda(z,t)$ is diagonal with components of order between $O(r^{-(m_j\pm\epsilon)})$.
\end{proposition}
\begin{proof}
Again we work in a small ball around the singularity $r=r(z,t)<\delta$ and this determines $\epsilon$.  As in the proof of Proposition~\ref{th:imf}, by replacing $g(z,t)$ by $g(z,t)p_i(z)$ we may assume $\gamma_i(z)=\lambda_i(z)$ is diagonal.  Using the notation $M$, $\Lambda$, $r_1$ and $r_2$ from the proof of Proposition~\ref{th:imf}, $g(t)=\exp(M(t))g(-\delta)$ for $M(t)=\Lambda(t)+M_0(t)$ where $M_0(t)$ is bounded and $\Lambda(t)$ is diagonal with asymptotic behaviour determined by the integral
\[\int^t_{-\delta}\frac{1}{r}dt<\int^t_{-\delta}\frac{1}{r_1}dt=\ln\frac{t+r_1}{-\delta+r_1(\delta)}
=\ln\frac{\delta+r_1(\delta)}{-t+r_1}
<-\ln r_1+B_1<-\ln r\]
and
\[\int^t_{-\delta}\frac{1}{r}dt>\int^t_{-\delta}\frac{1}{r_2}dt=\ln\frac{t+r_2}{-\delta+r_2(\delta)}
=\ln\frac{\delta+r_2(\delta)}{-t+r_2}
>-\ln 2r_2+B_2>-\ln r+B_2'\]
where we have used $0<-t<r_1<r<r_2$.  The last inequality in each uses $r_1\alpha>r>r_2/\alpha$ for some $\alpha>1$ which follows from 
\[r^2-r_1^2<(c_2-c_1)|z|^2<(c_2-c_1)\frac{r^2}{c_1}\quad\Rightarrow\quad (r_1/r)^2>2-c_2/c_1\]
\[r_2^2-r^2<(c_2-c_1)|z|^2<(c_2-c_1)\frac{r^2}{c_1}\quad\Rightarrow\quad (r_2/r)^2<c_2/c_1\]
and the fact that $c_1\approx c_2$.  These two inequalities can be represented heuristically by a straightforward infinitesimal inequality $d|t|<dr<d|t|+O(r)$.

Thus $\Lambda(z,t)$ has components lying between $(m_j\pm\epsilon)\ln r$.  The matrix $M$ is diagonalisable and its eigenspaces tend uniformly to the standard basis as $r\to 0$ with eigenvalues between $(m_j\pm\epsilon)\ln r$.  Thus $\exp(M)$ has eigenvalues of order $O(r^{m_j\pm\epsilon})$ as $r\to 0$ and eigenspaces tending uniformly to the standard basis.
\end{proof}

The spaces $\z(\lambda_1,y_1,...,\lambda_k,y_k)$ and $\y(\lambda_i,z_i)$ each have natural compactifications $\overline{\z(\lambda_1,y_1,...,\lambda_k,y_k)}$ and $\overline{\y(\lambda_i,z_i)}$ described as follows.  

The closure of the image $\y(\lambda_i,z_i)\hookrightarrow\Omega G$ gives a compactification $\overline{\y(\lambda_i,z_i)}$ with lower strata consisting of $\y(\lambda',z_i)$ for other conjugacy classes $\lambda':U(1)\to G$,  This induces a partial order on conjugacy classes of homomorphisms where $\lambda'\leq\lambda_i$ if $\y(\lambda',z_i)\subset\overline{\y(\lambda_i,z_i)}$.  (The partial order is independent of $z_i$.)

Any sequence of monopoles in $\z(\lambda_1,y_1,...,\lambda_k,y_k)$ has a subsequence that converges uniformly on compact subsets of $Y-\{y_1,...,y_k\}$ to a solution of the Bogomolny equation - "bubbling" occurs at the $y_i$ and nowhere else.  The boundary conditions are preserved so the the limit lies in $\z(\lambda_1',y_1,...,\lambda_k',y_k)$ for some $\lambda_i'$ and Kapustin and Witten showed that $\lambda_i'\leq\lambda_i$ using the ordering defined above and that this gives a compactification $\overline{\z(\lambda_1,y_1,...,\lambda_k,y_k)}$ with lower strata coinciding with the lower strata of $\overline{\y(\lambda_1,z_1)}\times...\times\overline{\y(\lambda_k,z_k)}$.
\begin{lemma}   \label{th:comp}
$\f$ extends continuously to \[\bar{\f}:\overline{\z(\lambda_1,y_1,...,\lambda_k,y_k)}\to\overline{\y(\lambda_1,z_1)}\times...\times\overline{\y(\lambda_k,z_k)}.\]
\end{lemma}
\begin{proof}
Solve the ODE (\ref{eq:scat}) that defines $\f$ along lines $(z,t)$ where $z$ is contained in a compact subset of $\Sigma-\{z_1,...,z_k\}$.  Since a sequence $(A_j,\Phi_j)$ of monopoles converges uniformly on compact sets, the corresponding sequence $g_j$ of solutions of (\ref{eq:scat}) satisfying $g_j|_{\Sigma\times 0}=I$ also converges uniformly, and hence so does its image under $\f$.
\end{proof}

\begin{cor}  \label{th:imc}
The image of $\f$ is closed.
\end{cor}
\begin{proof}
Take a component $\f_i$ of $\f=(\f_1,...,\f_k)$.  Let $\{\gamma_j\}\subset{\rm im}\f_i$ be a sequence in the image of $\f_i$ that converges $\{\gamma_j\}\to\gamma\in\y(\lambda_i,z_i)$.  Then there is a sequence $\{(A_j,\Phi_j)\}\subset\z(\lambda_1,y_1,...,\lambda_k,y_k)$ such that $\f_i(A_j,\Phi_j)=\gamma_j$.  This sequence posses\-ses a subsequence that converges in $\overline{\z(\lambda_1,y_1,...,\lambda_k,y_k)}$, $\{(A_{j_l},\Phi_{j_l})\}\to(A_{\infty},\Phi_{\infty})$.  Moreover, $\gamma_j=\f_i(A_j,\Phi_j)\to\f_i(A_{\infty},\Phi_{\infty})$ so $\f_i(A_{\infty},\Phi_{\infty})=\gamma_j$.  If the sequence had bubbled then $\f_i(A_{\infty},\Phi_{\infty})\in\y(\lambda',z_i)$ for $\lambda'<\lambda_i$ contradicting the assumption that $\gamma_i\in\y(\lambda_i,z_i)$.  Hence the sequence does not bubble so $(A_{\infty},\Phi_{\infty})\in\z(\lambda_i,y_i)$.  Thus the image of each $\f_i$ is closed so the image of $\f$ is closed.
\end{proof}
\noindent{\em Notes.}
\begin{enumerate}
\item[1.] Proposition~\ref{th:imf} was proven by Kronheimer when $G=SU(2)$ and the metric on $Y$ is flat (near $y_i$.)  The proof presented here is adapted from Kronheimer's proof.

\item[2.] Kapustin and Witten \cite{KWiEle} proved Proposition~\ref{th:imf} when $G=U(1)$ and hence for any diagonal monopole, by a topological argument.  Proposition~\ref{th:imf} confirms the intuition that the off-diagonal components of a monopole are negligible so it behaves as if it is diagonal.

\item[3.] The behaviour (\ref{eq:sing}) near singularities is a consequence of the ASD connection defined locally on $B^4$ above $B^3-y_i$.  Moreover it is equivalent to the existence of the ASD connection---given the behaviour (\ref{eq:sing}) one can show that the $L^2$ norm of the curvature of the ASD connection upstairs is finite.  Then apply Uhlenbeck's result over a 4-dimensional base to get $C^0$ bounds on the curvature from $L^2$ bounds on the curvature (and also from these bounds finds good gauges in which $A$ itself obeys $C^0$ bounds.)  This was proven by Kronheimer \cite{KroMon} in the Euclidean case and by Pauly \cite{PauMon} in the general case when the metric on the 4-manifold upstairs is not necessarily smooth, (for $G=SU(2)$ although the proof generalises immediately.)  In particular,
the conditions (\ref{eq:sing}) control the growth of the entire monopole.  In the abelian case it is straightforward that (\ref{eq:sing}) determines the behaviour of the monopole near $y_i$, or more precisely there exists a gauge in which the growth of $A$ is controlled.  It is similarly straightforward from (\ref{eq:sing}) to determine the behaviour of the diagonal part of $A$ (\ie the part parallel to $\Phi$), while the off-diagonal components of $A$ require Uhlenbeck's deep results.

\item[4.]  Traditionally, non-abelian magnetic monopoles have been studied only over non-compact $Y$.   Monopoles over compact $Y$ reduce to the flat case in which the moduli space is essentially a discrete set.  This is because the Higgs field satisfies a maximum principle
\[ \quad\quad\quad\quad\Delta|\Phi|^2=2*d*\langle d^A\Phi,\Phi\rangle=
2*d\langle F_A,\Phi\rangle=2*\langle F_A,d^A\Phi\rangle=2|F_A|^2\geq 0.\]
The introduction of singularities enables non-trivial solutions over compact $Y$.   The monopole boundary value problem generalises the abelian case of harmonic function boundary value problems.

\item[5.]  A consequence of (\ref{eq:sing}) is
\[ F_A=O\left(\frac{1}{r^2}\right)\]
so in particular the curvature on the 3-manifold is not square-integrable.  Thus its integral is not a topological invariant, in contrast to the non-singular case where $||F_A||_2$ is finite and a topological invariant that distinguishes components of the moduli space.

\item[6.]  Singular $U(1)$ monopoles on $Y^3$ have three interpretations.  Firstly as Dirac monopoles, secondly as Green's functions and thirdly as metrics on a four-manifold with a $U(1)$ action and quotient $Y^3$.  The existence of singular $U(1)$ monopoles comes from the existence of Green's functions which is deep---see for example \cite{TayPar}.  Via a homomorphism $U(1)\to G$, singular $U(1)$ monopoles are used in the existence proof for general $G$ singular monopoles.  Furthermore, singular $U(1)$ monopoles are required for two applications of the maximum principle---to show that $\f$ has open image, and that it is one-to-one onto its image---in order to remove singularities of subharmonic functions.  It is worth including the following easy fact.
\begin{lemma}
A Green's function on $\Sigma\times I-\{y_1,...,y_n\}$ uniquely determines a Dirac monopole satisfying the boundary conditions (\ref{eq:bou}).
\end{lemma}
\begin{proof}
Denote the Green's function by $V$, so $\Phi=iV$.  If two $U(1)$ connections $A_1$ and $A_2$ satisfy $dA_i=*d\Phi$, their difference $a=A_1-A_2$ is a closed 1-form.  In particular, the integral of $a$ around closed loops is a well-defined function on homotopy classes.  The integral of $a$ around any closed loop on $\Sigma\times0$ is zero since $A_i|_{\Sigma\times0}=0$.  The inclusion $\Sigma\times0\hookrightarrow\Sigma\times I-\{y_1,...,y_n\}$ induces an isomorphism on the fundamental group hence the integral of $a$ around any closed loop is zero and $a$ is exact.  Thus the connections are gauge equivalent under the gauge transformation $\exp{f}$ where $a=df$.  
\end{proof}

\end{enumerate}

\section{Linearisation}   \label{sec:lin}

Kapustin and Witten \cite{KWiEle} showed that the tangent spaces of $\z(\lambda_1,y_1,...,\lambda_k,y_k)$ and $\y(\lambda_1,z_1)\times\y(\lambda_k,z_k)$ have the same dimension.  This follows from an index calculation given in \cite{KWiEle} for $Y=\Sigma\times I$ that generalises a result of Pauly \cite{PauMon} when the 3-manifold is closed and $G=SU(2)$.  The tangent space $T_{(A,\Phi)}\z(\lambda,y_0)$ is given by $H^1$ of the complex
\[ 0\to\Omega^0(\ad(E))\stackrel{d_{A,\Phi}}{\longrightarrow}\Omega^1(\ad(E))\oplus\Omega^0(\ad(E))\stackrel{L_{A,\Phi}}{\longrightarrow}\Omega^2(\ad(E))\to 0\]
where $d_{A,\Phi}\xi=(d^A\xi,[\Phi,\xi])$ is an infinitesimal gauge transformation of $(A,\Phi)$ and
\begin{equation}  \label{eq:boglin} 
L_{A,\Phi}(a,\phi):=d^Aa+*[\Phi,a]-*d^A\phi=0
\end{equation}
is the linearisation of the Bogomolny equation.  The boundary conditions require the infinitesimal gauge transformations $\xi\in\Omega^0(\ad(E))$ and the 1-form $a\in\Omega^1(\ad(E))$ to vanish on the boundary component $\Sigma\times0$ and the 0-form $\phi\in\Omega^0(\ad(E))$ to vanish on the boundary component $\Sigma\times0$.  

Using the product structure $\Sigma\times I$ the linearised Bogomolny equation (\ref{eq:boglin}) can be expressed as the linearisation of (\ref{eq:bogprod})
\begin{equation*}  \label{eq:boglinprod}
\begin{array}{rcl}
\partial_{\bar{z}}^A(a_t-i\phi)&=&(\partial_t^A-i\Phi)a_{\bar{z}}\\
\partial_{\bar{z}}^Aa_z-\partial_z^Aa_{\bar{z}}&=&\rho(z,\bar{z})\left((\partial_t^A+i\Phi)(a_t-i\phi)-(\partial_t^A-i\Phi)(a_t+i\phi)\right)
\end{array}\tag{6a}
\end{equation*}
where $a_{\bar{z}}$ is the $d\bar{z}$ component of $a$.  

When $G=SU(2)$, Pauly \cite{PauMon} set up a Banach manifold that contains the singular monopoles and showed that $L_{A,\Phi}\oplus d_{A,\Phi}^*$ is Fredholm when the base manifold is closed with index $4m$.  Pauly's argument uses excision properties of the index allowing one to localise the calculation near the singularity and pass to the singular $U(1)$ bundle over a neighbourhood of the singularity.  Given any number of singularities, of weights $m_i$, the argument shows that the index is additive, yielding $4\sum_im_i$.  Kapustin and Witten \cite{KWiEle} observed that the excision argument allows the calculation to generalise to manifolds with boundary.  They showed that in the absence of singularities, a monopole satisfying the boundary conditions (\ref{eq:bou}) is trivial (this also follows from the maximum principle argument in Section~\ref{sec:uniq}) and the index in this case is zero.  When a singularity is added the index becomes $0+4m$.  Furthermore they showed it is straightforward to generalise to any $G$ since the calculation reduces to a known equivariant calculation over $S^4$.

The boundary conditions ensure that the map $\xi\mapsto(d^A\xi,[\Phi,\xi])$ is injective and $L_{A,\Phi}$ is surjective.  
Hence the tangent space $T_{(A,\Phi)}\z(\lambda_1,y_1,...,\lambda_k,y_k)$ is defined from the complex as the quotient of the kernel of $L_{A,\Phi}$ by the gauge directions, with dimension equal to the index, and this coincides with the real dimension of $\y(\lambda_1,z_1)\times...\times\y(\lambda_k,z_k)$.

The derivative 
\[ D\f_{(A,\Phi)}:T_{(A,\Phi)}\z(\lambda_1,y_1,...,\lambda_k,y_k)\to \bigoplus_{i=1}^k T_{\gamma_i}\y(\lambda_i,z_i)\] 
for $\gamma_i=\f_i(A,\Phi)$ is  is determined by the linearisation of (\ref{eq:scat})
\begin{equation}  \label{eq:scatlin}
 (\partial_t^A-i\Phi)\eta=-(a_t-i\phi)
 \end{equation}
for $\eta$ a locally defined ${\bf g}^c$ valued function.  (So $A_t-i\Phi$ acts by Lie bracket on $\eta$.)  The $i$th component of $D\f_{(A,\Phi)}(a,\phi)$ is the derivative of $\f_i$, notated by $D\f_i(a,\phi)$ (suppress $(A,\Phi)$), and determined by local solutions of (\ref{eq:scatlin}) near $y_i$.  The following lemma defines $D\f_i(a,\phi)$ from $\eta$ and $g(z,t)$, the solution of (\ref{eq:scat}) satisfying $g|_{\Sigma\times 0}=I$.
Let $U_i\subset\Sigma$ be a disk neighbourhood of $z_i$.
\begin{lemma}  \label{th:defdf}
The two solutions of (\ref{eq:scatlin}), $\eta:\Sigma\times I-\displaystyle\bigcup_i\{z_i\}\times(t_i,1]\to{\bf g}^c$ such that $\eta(0)=0$ and $\eta_i:U_i\times I-\{z_i\}\times[0,t_i)\to{\bf g}^c$  satisfying $\partial_{\bar{z}}^A\eta=\partial_{\bar{z}}^A\eta_i$ are related by
\begin{equation}  \label{eq:etas} 
\eta_i(z,t)=\eta(z,t)+g(z,t)\xi_i(z)g(z,t)^{-1}
\end{equation}
for $\xi_i(z)\in L^+{\bf g}^c$.   This defines the derivative at $(A,\Phi)$
\[
D\f_i(a,\phi)=\xi_i(z)\gamma_i(z)\in T_{\gamma_i}\Omega G.
\]
\end{lemma}
\begin{proof}
 The fundamental solution of the homogeneous equation associated to (\ref{eq:scatlin}) is $g(z,t)\xi(z)g(z,t)^{-1}$ where $g(z,t)$ satisfies (\ref{eq:scat}) and $\xi(z)$ is independent of $t$.  Thus since $\eta$ and $\eta_i$ are two solutions of (\ref{eq:scatlin}) they satisfy (\ref{eq:etas}) for some $\xi_i:U-\{z_i\}\to{\bf g^c}$.  It remains to prove that $\xi_i$ is holomorphic and extends over $U$.

Differentiate (\ref{eq:scatlin}) and use $[\partial_{\bar{z}}^A,\partial_t^A-i\Phi]=0$ from (\ref{eq:bogprod}), to get
\[ (\partial_t^A-i\Phi)\partial_{\bar{z}}^A\eta=\partial_{\bar{z}}^A(\partial_t^A-i\Phi)\eta=-\partial_{\bar{z}}^A(a_t-i\phi)\]
which combined with the first equation of (\ref{eq:boglinprod}) yields
\[ (\partial_t^A-i\Phi)\partial_{\bar{z}}^A\eta=-(\partial_t^A-i\Phi)a_{\bar{z}}.\]
Thus
\begin{equation}   \label{eq:az}
\partial_{\bar{z}}^A\eta=-a_{\bar{z}}.
\end{equation}
since $(\partial_t^A-i\Phi)(\partial_{\bar{z}}^A\eta+a_{\bar{z}})=0$ and at $t=0$, $\partial_{\bar{z}}^A\eta=0=-a_{\bar{z}}$.

Differentiate (\ref{eq:etas}) by $\partial_{\bar{z}}^A$ to get
\[ \partial_{\bar{z}}^A\eta_i(z,t)=-a_{\bar{z}}+g(z,t)\partial_{\bar{z}}\xi_i(z)g(z,t)^{-1}\]
and since $\partial_{\bar{z}}^A\eta_i(z,t)=\partial_{\bar{z}}^A\eta(z,t)=-a_{\bar{z}}$ (on the overlap of their domains) this implies $\partial_{\bar{z}}\xi_i(z)=0$.

The Lie algebra valued functions $\eta$ and $\xi$ arise as variations of $g$ and $\gamma$
\[\delta g=\eta g,\quad\delta\gamma_i=\xi_i\gamma_i.\]
In particular, $T_{\gamma}\y(\lambda,z_i)=L^+{\bf g}^c\cdot\gamma$ so if $(a,\phi)\in T_{(A,\Phi})\z(\lambda_1,y_1,...,\lambda_k,y_k)$ then $\xi_i\gamma_i\in T_{\gamma}\y(\lambda,z_i)$ and we can choose $\xi_i\in L^+{\bf g}^c$.
\end{proof}

\begin{proposition}   \label{th:openim}
The map $\f:\z(\lambda_1,y_1,...,\lambda_k,y_k)\to\y(\lambda_1,z_1)\times...\times\y(\lambda_k,z_k)$ has open image.
\end{proposition}
\begin{proof}
We need to prove that $D\f_{(A,\Phi)}$  is an isomorphism for each $(A,\Phi)$ so that the inverse function theorem applies to give the result.  As described above, the index theorem together with regularity of $L_{A,\Phi}$ gives the fact that the dimensions of the spaces coincide.  Hence, it is sufficient to prove that $D\f_{(A,\Phi)}$ is injective which will follow from the injectivity of $D\f_i$ for each $i=1,...,k$.  The proof of this is an application of the maximum principle.

\begin{lemma}  \label{th:ellin}
\[\left(\partial_{\bar{z}}^A\partial_z^A+\rho(z,\bar{z})(\partial_t^A-i\Phi)(\partial_t^A+i\Phi)\right)(\eta+\eta^*)=0\]
\end{lemma}
\begin{proof}
Into the second equation of (\ref{eq:boglinprod}), substitute (\ref{eq:scatlin}) and (\ref{eq:az}) together with their adjoints $(\partial_t^A+i\Phi)\eta^*=a_t+i\phi$ and $\partial_z^A\eta^*=a_z$ to get
\[ \partial_{\bar{z}}^A\partial_z^A\eta^*+\partial_z^A\partial_{\bar{z}}^A\eta=-\rho(z,\bar{z})\left((\partial_t^A+i\Phi)(\partial_t^A-i\Phi)\eta+(\partial_t^A-i\Phi)(\partial_t^A+i\Phi)\eta^*\right)\]
and thus
\[\left(\partial_{\bar{z}}^A\partial_z^A+\rho(\partial_t^A-i\Phi)(\partial_t^A+i\Phi)\right)(\eta+\eta^*)=\left[[\partial_{\bar{z}}^A,\partial_z^A]+\rho[\partial_t^A-i\Phi,\partial_t^A+i\Phi],\eta\right].\]
The zero order operator on the right hand side vanishes by the Bogomolony equation and the lemma is proven.
\end{proof}
\begin{lemma}  \label{th:subh}
$\Delta|\eta+\eta^*|^2\geq 0$
\end{lemma}
\begin{proof}
We have
$|\eta+\eta^*|^2=\langle\eta+\eta^*,\eta+\eta^*\rangle$ and $\Delta=(1/\rho)\partial_{\bar{z}}\partial_z+\partial^2_t$ so
\begin{eqnarray*}
\Delta|\eta+\eta^*|^2&=&\frac{1}{\rho}\partial_{\bar{z}}\partial_z\langle\eta+\eta^*,\eta+\eta^*\rangle+\partial^2_t\langle\eta+\eta^*,\eta+\eta^*\rangle\\
&=&2Re\langle\frac{1}{\rho}\partial_{\bar{z}}^A\partial_z^A(\eta+\eta^*),\eta+\eta^*\rangle+\frac{2}{\rho}\left|\partial_z^A(\eta+\eta^*)\right|^2\\
&+&\hspace{-.5mm}2Re\langle(\partial_t^A-i\Phi)(\partial_t^A+i\Phi)(\eta+\eta^*),\eta+\eta^*\rangle+2\left|(\partial_t^A+i\Phi)(\eta+\eta^*)\right|^2\\
&\geq&2Re\langle\left(\frac{1}{\rho}\partial_{\bar{z}}^A\partial_z^A+(\partial_t^A-i\Phi)(\partial_t^A+i\Phi)\right)(\eta+\eta^*),\eta+\eta^*\rangle=0
\end{eqnarray*}
where the final step uses Lemma~\ref{th:ellin}.
\end{proof}

We would like to use Lemma~\ref{th:subh} to apply the maximum principle to $|\eta+\eta^*|^2$, but there are two further critical issues we must address.  Firstly a solution of (\ref{eq:scatlin}) along $z=z_i$ does not extend past $t_i$ so in general $\eta$ and hence $|\eta+\eta^*|^2$ is not globally defined on $Y=\Sigma\times I-\{y_1,...,y_k\}$.  Secondly, to extract something meaningful from the maximum principle, we need to know the growth of $|\eta+\eta^*|^2$ near each $y_i$.

We deal with these two issues in the following two lemmas.
\begin{lemma}  \label{th:etaextends}
For $\eta$ solving (\ref{eq:scatlin}) and satisfying $\eta|_{\Sigma\times0}=0$,
\[ D\f_{(A,\Phi)}(a,\phi)=0\quad{\rm implies}\quad\eta:Y-\{y_1,...,y_k\}\to{\bf g}^c.\]
\end{lemma}
\begin{proof}
The function $\eta$ is defined on $\Sigma\times I-\cup_iz_i\times(t_i,1]$.  The kernel of $D\f_{(A,\Phi)}$ is the intersection of the kernels of $D\f_i$ at $(A,\Phi)$.  We will show that for each $i$, if $D\f_i(a,\phi)=0$ then $\eta$ extends to $z_i\times(t_i,1]$.  If this is true for all $i=1,...,k$ then $\eta$ extends to all of  $\Sigma\times I-\{y_1,...,y_k\}$.

From Lemma~\ref{th:defdf}, $\xi_i(z)\in L^+{\bf g^c}$ defined by (\ref{eq:etas}) represents the zero vector if 
\[\gamma_i(z)^{-1}\xi_i(z)\gamma_i(z)=\xi_+\in L^+{\bf g^c}.\] 
Then (\ref{eq:etas}) becomes $\eta_i(z,t)=\eta(z,t)+g(z,t)\gamma_i(z)\xi_+(z)\gamma_i(z)^{-1}g(z,t)^{-1}$.  By definition $g(z,t)\gamma_i(z)$ is defined on $U_i\times(t_i,1]$, hence $g(z,t)\gamma_i(z)\xi_+(z)\gamma_i(z)^{-1}g(z,t)^{-1}$ is defined on $U_i\times(t_i,1]$.  Since $\eta_i$ is also defined on $U_i\times(t_i,1]$, (\ref{eq:etas}) implies $\eta$ is defined on $U_i\times(t_i,1]$ and hence extends.  This is is true for each $i=1,...,k$, so $\eta$ is well-defined globally on $\Sigma\times I-\{y_1,...,y_k\}$, if $(a,\phi)\in\ker D\f_{(A,\Phi)}$.  
\end{proof}

\begin{lemma}
For $\eta$ solving (\ref{eq:scatlin}) and satisfying $\eta|_{\Sigma\times0}=0$,
\[ D\f_{(A,\Phi)}(a,\phi)=0\quad{\rm implies}\quad|\eta|=O(r^{-\epsilon}){\rm\ near\ each\ }y_i,\quad{\rm for\ }\epsilon>0.\]
\end{lemma}
\begin{proof}
The asymptotic behaviour of $\eta$ depends only on small balls around each $y_i$.  As in the proofs of Propositions~\ref{th:imf} and \ref{th:imfc} we choose small balls of radius $\delta$ around each $y_i$ which determines $\epsilon$.  Equation (\ref{eq:scatlin}) is linear in $\eta$ so we can study the growth of each component of $\eta$ separately.  The idea is to separate $\eta$ into components of two types which have small growth for different reasons.  Components of one type have growth $o(r^{-\epsilon})$ simply due to the fact that $\eta$ extends to $Y-\{y_1,...,y_k\}$ by Lemma~\ref{th:etaextends}.  Components of the other type have growth $o(r^{-\epsilon})$ due to the ODE (\ref{eq:scatlin}) together with the fact $a_t-i\phi=O(1/r)$ near each singularity which follows, as in the proof of (\ref{eq:sing}), from the continuity of a Lie algebra valued 1-form upstairs.   

Since $g(z,t)\sim\exp\Lambda(z,t)\cdot p(z)^{-1}$ as $r\to 0$ where $\exp\Lambda(z,t)$ is diagonal with components of order between $O(r^{-(m_j\pm\epsilon)})$ as $r\to 0$, we can choose a basis $\{b_{jk}(z)\}$ of ${\bf g}^c$ obtained by conjugating the standard basis $e_{jk}$ by $p(z)$ with the property that the component $g(z,t)b_{jk}(z)g(z,t)^{-1}$ behaves like $r^{m_k-m_j\pm 2\epsilon}$.  In particular, this component grows when $m_j>m_k$.  Put 
\begin{equation}   \label{eq:decomp}
\eta(z,t)=\sum\eta_{jk}(t)g(z,t)b_{jk}(z)g(z,t)^{-1}.
\end{equation}  
Then $\eta g\gamma_i=g\gamma_i\sum\eta_{jk}(t)z^{m_k-m_j}e_{jk}$.  But $g\gamma_i$  and $\eta$ are defined on $t>t_i$ by definition of $\gamma_i$ and by Lemma~\ref{th:etaextends} since $D\f_{(A,\Phi)}(a,\phi)=0$.  Thus, when $m_j>m_k$, $\eta_{jk}(t)$ must decay at least as fast as $t^{m_j-m_k}$ so that $\eta_{jk}(t)z^{m_k-m_j}$ is defined.  (We have changed coordinates $t\mapsto t-t_i$.)  In particular, the component $\eta_{jk}(t)g(z,t)b_{jk}(z)g(z,t)^{-1}=O(r^{\pm 2\epsilon})$ when $m_k<m_j$.

Using (\ref{eq:decomp})
\[0=(\partial_t^A-i\Phi)\eta+(a_t-i\phi)=\sum\frac{d\eta_{jk}}{dt}g(z,t)b_{jk}(z)g(z,t)^{-1}+O(1/r).\]
When $m_j\leq m_k$, $gb_{jk}g^{-1}=O(r^{m_j-m_k\pm 2\epsilon})$ so $d\eta_{jk}/dt=O(r^{m_j-m_k-1\pm 2\epsilon})$ which grows as $r\to 0$.  By integration $\eta_{jk}=O(r^{m_j-m_k\pm 2\epsilon})$ and thus the component $\eta_{jk}g(z,t)b_{jk}(z)g(z,t)^{-1}=O(r^{\pm 4\epsilon})$ as required.  Notice that the second argument requires $m_j\leq m_k$ since if $m_j>m_k$ then by integration we can only conclude that $\eta_{jk}=O(1)$ whereas better decay is needed.
\end{proof}
The function $|\eta+\eta^*|^2$ is well-defined on $Y-\{y_1,...,y_k\}$ and near each $y_i$ since $|\eta|=o(r^{-\epsilon})$ then $|\eta+\eta^*|^2=o(r^{-2\epsilon})=o(1/r)$ for $\epsilon$ small enough.  Thus $|\eta+\eta^*|^2$ is subharmonic and satisfies $|\eta+\eta^*|^2|_{\Sigma\times 0}=0$.  At $\Sigma\times 1$, $\Phi=0=\phi$ so $\partial^A_t\eta=-a_t$ implies $\partial_t|\eta+\eta^*|^2=0$.  Take a positive harmonic function $V$ satisfying $V|_{\Sigma\times 0}=0$, $\partial_tV|_{\Sigma\times 1}=0$ and $V=O(1/r)$ near each $y_i$.  It is given by a sum of Green's functions for the mixed Dirichlet/Neumann harmonic function boundary value problem.  The function $V$ grows faster than $|\eta+\eta^*|^2$ near $r=0$ so $|\eta+\eta^*|^2-cV\leq 0$ near $r=0$ for all $c>0$.  By the maximum principle, $|\eta+\eta^*|^2-cV\leq 0$ for all $c>0$ thus $|\eta+\eta^*|^2\leq 0$, and hence $|\eta+\eta^*|^2=0$, or equivalently $\eta=-\eta^*$ lives in ${\bf g}$.  But $\eta$ determines $(a,\phi)$, by (\ref{eq:scatlin}) and (\ref{eq:az}) and if $\eta\in{\bf g}$ then $(a,\phi)=(-d^A\eta,-[\Phi,\eta])$ which corresponds to an infinitesimal gauge transformation of $(A,\Phi)$.  Directions along the gauge orbit are trivial in the tangent space of the moduli space so $D\f$ is injective and the proposition is proven.
\end{proof}

\begin{cor}   \label{th:surj}
The map $\f$ is surjective.
\end{cor}
\begin{proof}
The proof uses the method of continuity.  The existence of Green's functions and a homomorphism $U(1)\to G$ of weight $\lambda$ proves that $\z(\lambda_1,y_1,...,\lambda_k,y_k)$ is non-empty.  By Proposition~\ref{th:openim} the image of $\f$ is open.  The image of $\f$ is closed by Corollary~\ref{th:imc}, hence it is open, closed and non-empty so by the connectedness of the target space it is all of $\y(\lambda_1,z_1)\times...\times\y(\lambda_k,z_k)$.
\end{proof}

\section{Uniqueness}   \label{sec:uniq}

Consider the homogeneous space $G\backslash G^c$ given by the quotient of $G^c$ by a left action of $G$.  When $G=GL(n,\bc)$ this is is the space of Hermitian matrices.  Define a Hermitian metric $H$ on a trivial $G^c$ bundle $E$ over $Y-\{y_1,...,y_k\}$ to be a section of the trivial $G\backslash G^c$ bundle.  Note that $H$ may require more than one chart to be defined since it takes values in a quotient space.  With respect to a product structure $Y=\Sigma\times I$, a Hermitian metric $H$ on $E$ determines a monopole on $E$ as follows:
\[ \partial^A_t-i\Phi=\partial_t,\quad \partial^A_{\bar{z}}=\partial_{\bar{z}},\quad
\partial^A_t+i\Phi=\partial_t+H^{-1}\partial_tH,\quad \partial^A_z=\partial_z+H^{-1}\partial_zH\]
which is the monopole expressed in a gauge in which $\partial^A_t-i\Phi$ and $\partial^A_{\bar{z}}$ are trivial.  Hence to determine $H$ from $(A,\Phi)$, we express the monopole in the gauge $g$ satisfying (\ref{eq:scat}) and define $H$ to be the image of $g$ in the quotient $G\backslash G^c$ bundle, or equivalently $H=g^*g$.  In this gauge, the Bogomolny equation is
\[ \partial_{\bar{z}}(H^{-1}\partial_zH)+\rho(z,\bar{z})\partial_t(H^{-1}\partial_tH)=0.\]
In this case $H$ requires more than one chart to be defined since $g$ is not globally defined as it does not extend to $z_i\times (t_i,1]$ for $i=1,...,k$.  Recall that $g\gamma_i$ is well-defined on a neighbourhood of $z_i\times (t_i,1]$.  Thus $\bar{\gamma}_i^TH\gamma_i$ is also defined on a neighbourhood of $z_i\times (t_i,1]$ and represents $H$ there. (We have abused notation and used $H$ for a Hermitian metric and one of its local representatives.)  Following Donaldson \cite{DonAnt} we use the metric-like function on pairs of Hermitian metrics 
\[\sigma(H_1,H_2)=\tr H_1^{-1}H_2+\tr H_2^{-1}H_1-2n\] 
where $n$ is the rank of $G$. It has the property that $\sigma(H_1,H_2)$ vanishes if and only if $H_1=H_2$.

\begin{lemma}
The complex gauge transformation $h=H_1^{-1}H_2$ satisfies the equation 
\begin{equation}   \label{eq:endo} 
\partial_{\bar{z}}^{A_1}(h^{-1}\partial_z^{A_1}h)+\rho(z,\bar{z})(\partial_t^{A_1}-i\Phi)(h^{-1}(\partial_t^{A_1}+i\Phi_1)h)=0.
\end{equation}
\end{lemma}
\begin{proof}
The expression $h=H_1^{-1}H_2$ makes sense because we identify the gauges in which $\partial_t^{A_1}-i\Phi_1=\partial_t=\partial_t^{A_2}-i\Phi_2$.  We have 
\begin{eqnarray*}
\partial_z^{A_1}+h^{-1}\partial_z^{A_1}h&=&\partial_z^{A_1}+h^{-1}\partial_zh+h^{-1}[H_1^{-1}\partial_zH_1,h]\\
&=&\partial_z^{A_2}
\end{eqnarray*}
hence
\[ h^{-1}\partial_z^{A_1}h=-H_1^{-1}\partial_zH_1+H_2^{-1}\partial_zH_2\]
and
\[\partial_{\bar{z}}^{A_1}(h^{-1}\partial_z^{A_1}h)=-\partial_{\bar{z}}(H_1^{-1}\partial_zH_1)+\partial_{\bar{z}}(H_2^{-1}\partial_zH_2)\]
where we have used $\partial_{\bar{z}}^{A_1}=\partial_{\bar{z}}$ in the last step to get a gauge independent expression.  Also
\begin{eqnarray*}
\partial_t^{A_1}+i\Phi_1+h^{-1}(\partial_t^{A_1}+i\Phi_1)h&=&\partial_t^{A_1}+i\Phi_1+h^{-1}\partial_th+h^{-1}[H_1^{-1}\partial_tH_1,h]\\
&=&\partial_t^{A_2}+i\Phi_2
\end{eqnarray*}
hence
\[h^{-1}(\partial_t^{A_1}+i\Phi_1)h=-H_1^{-1}\partial_tH_1+H_2^{-1}\partial_tH_2\]
and
\[(\partial_t^{A_1}-i\Phi)(h^{-1}(\partial_t^{A_1}+i\Phi_1)h)=-\partial_t(H_1^{-1}\partial_tH_1)+\partial_t(H_2^{-1}\partial_tH_2)\]
where we have used $\partial_t^{A_1}-i\Phi_1=\partial_t$ in the last step.  

Thus, if $(A_1,\Phi_1)$ and $(A_2,\Phi_2)$ are two monopoles, adding these two equations and using the vanishing
\[ \partial_{\bar{z}}(H_j^{-1}\partial_zH_j)+\rho(z,\bar{z})\partial_t(H_j^{-1}\partial_tH_j)=0,\quad j=1,2\] 
then $h$ satisfies the gauge independent equation 
\[ \partial_{\bar{z}}^{A_1}(h^{-1}\partial_z^{A_1}h)+\rho(z,\bar{z})(\partial_t^{A_1}-i\Phi)(h^{-1}(\partial_t^{A_1}+i\Phi_1)h)=0.\]
\end{proof}

\begin{lemma}
With $h = H_1^{-1}H_2$ as above coming from two monopoles, 
\[ \Delta(tr(h)) \geq 0\] 
on $\Sigma\times I$. 
\end{lemma}
\begin{proof}
Differentiate $tr(h)$ to get
\[
\partial_ztr(h)=tr(\partial_z^Ah)=tr(h\cdot h^{-1}\partial_z^Ah)\]
so
\begin{eqnarray*}
\partial_{\bar{z}}\partial_ztr(h)&=&
tr\left((\partial_{\bar{z}}^Ah)h^{-1}\partial_z^Ah\right)+tr\left(h\partial_{\bar{z}}^A(h^{-1}\partial_z^Ah)\right)\\
&\geq& tr\left(h\partial_{\bar{z}}^A(h^{-1}\partial_z^Ah)\right)
\end{eqnarray*}
where the inequality uses $tr\left((\partial_{\bar{z}}^Ah)h^{-1}\partial_z^Ah\right)\geq 0$ since $\partial_z^Ah=(\partial_{\bar{z}}^Ah)^*$.  Similarly,
\[\partial_ttr(h)=tr((\partial_t^A+i\Phi)h)=tr(h\cdot h^{-1}(\partial_t^A+i\Phi)h)\]
so
\begin{eqnarray*}
\partial_t^2tr(h)&=&
\hspace{-.5mm}tr(((\partial_t^A-i\Phi)h) h^{-1}(\partial_t^A+i\Phi)h)+
tr(h(\partial_t^A-i\Phi)(h^{-1}(\partial_t^A+i\Phi)h))\\
&\geq& tr(h(\partial_t^A-i\Phi)(h^{-1}(\partial_t^A+i\Phi)h))
\end{eqnarray*}
the inequality again using the fact that $(\partial_t^A-i\Phi)h=((\partial_t^A+i\Phi)h)^*$.  Thus
\begin{eqnarray*}
\rho(z,\bar{z})\Delta(tr(h))&=&\partial_{\bar{z}}\partial_ztr(h)+\rho(z,\bar{z})\partial_t\partial_ttr(h)\\
&\geq&tr [h\partial_{\bar{z}}^A(h^{-1}\partial_z^Ah)+h\rho(z,\bar{z})(\partial_t^A-i\Phi)(h^{-1}(\partial_t^A+i\Phi)h)]=0
\end{eqnarray*}
by (\ref{eq:endo}).
\end{proof}

\begin{lemma}  \label{th:sigo1}
If $\f(A_1,\Phi_1)=\f(A_2,\Phi_2)$ then the $H_j$ corresponding to $(A_j,\Phi_j)$ satisfy $\sigma(H_1,H_2)=o(1/r)$ near each singularity $y_i$.
\end{lemma}
\begin{proof}
Consider each component $\f_i$ of $\f=(\f_1,...,\f_k)$ separately.  By Proposition~\ref{th:imfc}, $\f_i(A_1,\Phi_1)=\f_i(A_2,\Phi_2)=\gamma_i(z)=p_i(z)\lambda_i(z)$ implies that the solutions $g_1(z,t)$ and $g_2(z,t)$ of (\ref{eq:scat}) are asymptotic to
$\exp\Lambda_1(z,t)\cdot p_i(z)^{-1}$, respectively $\exp\Lambda_2(z,t)\cdot p_i(z)^{-1}$, where the diagonals $\exp\Lambda_1$ and $\exp\Lambda_2$ have components of order between $O(r^{-(m_j\pm\epsilon)})$.  Thus $||g_2g_1^{-1}||=O(r^{-2\epsilon})$.
Now $H_j=g_j^*g_j$ so $tr(H_1^{-1}H_2)=tr(g_1^{-1}(g^*_1)^{-1}g_2^*g_2)=tr(g_2g_1^{-1}(g^*_1)^{-1}g_2^*)=||g_2g_1^{-1}||^2=O(r^{-4\epsilon})=o(1/r)$.  Similarly for $tr(H_2^{-1}H_1)$, thus $\sigma(H_1,H_2)=o(1/r)$ near each singularity $y_i$.
\end{proof}

\begin{proposition}   \label{th:1to1}
The map $\f:\z(\lambda_1,y_1,...,\lambda_k,y_k)\to\y(\lambda_1,z_1)\times...\times\y(\lambda_k,z_k)$ is 1-1 onto its image.
\end{proposition}
\begin{proof}
Again consider each component $\f_i$ separately.  Suppose $\f_i(A_1,\Phi_1)=\gamma_i=\f_i(A_2,\Phi_2)$.  As for the linear case, there are two critical issues we must address in order to apply the maximum principle.  Firstly a solution of (\ref{eq:scatlin}) along $z=z_i$ does not extend past $t_i$ so in general $\sigma(H_1,H_2)$ is not globally defined on $Y=\Sigma\times I-\{y_1,...,y_k\}$.  Secondly, to extract something meaningful from the maximum principle, we need to know the growth of $\sigma(H_1,H_2)$ near each $y_i$.

We resolve the first issue using $H_j\mapsto\bar{\gamma}_i^TH_j\gamma_i$ for $j=1,2$ near $y_i$ so $H_1^{-1}H_2\mapsto\gamma_i^{-1}H_1^{-1}H_2\gamma_i$.  Thus $\sigma(H_1,H_2)=\tr H_1^{-1}H_2+H_2^{-1}H_1-4$ is  a globally well-defined function on $Y-\{y_1,...,y_k\}$.   For the second issue, we have $\sigma(H_1,H_2)=o(1/r)$ by Lemma~\ref{th:sigo1}.

Now $\Delta\sigma\geq 0$ so $\sigma$ is subharmonic.  As in the linear case, $\sigma|_{\Sigma\times 0}=0$ and $\partial_t\sigma|_{\Sigma\times 1}=0$ so by comparison with a positive harmonic function $V$ satisfying the same boundary conditions and $O(1/r)$ near the singularities, $\sigma-cV\leq 0$ for all $c>0$ since $O(1/r)>o(1/r)$ thus $\sigma\leq 0$.  Hence $\sigma(H_1,H_2)=0$ and  $H_1=H_2$.
\end{proof}

\section{Dimensional reduction of four-manifolds.}

Singular monopoles over a 3-manifold $Y$ can be described locally as $U(1)$ invariant instantons on a local Riemannian 4-manifold.  In this section we construct a {\em global} Riemannian 4-manifold $X$ equipped with a circle action so that $X/U(1)\cong\Sigma\times I$.  The action is free away from fixed points $\hat{p}_j$ lying above the singular points $y_j\in Y$ and gives the Hopf fibration in a neighbourhood above each singular point $y_j$.  Singular monopoles on the trivial bundle $E$ over the product $\Sigma\times I$ pull back to $U(1)$ invariant instantons on a bundle $\tilde{E}$ over $X$.  The weight of the action of $U(1)$ on $\tilde{E}|_{\hat{p}_j}$ is the weight of the singularity of the monopole at $y_j=(z_j,t_j)$.  This point of view leads to a second proof of uniqueness and existence of singular monopoles.  Note that there is a topological obstruction to the existence of $X$ when $Y$ has no boundary, \cite{PauMon}.  There is a choice of topological type of the manifold $X$ and for each such choice, there are many Riemannian metrics on $X$ with quotient the given product metric on $\Sigma\times I$, so we choose one.  The Riemannian manifold $X$ will be complex and K\"ahler and for any Hecke modification over $Y=\Sigma\times I$ one can construct a holomorphic bundle over $X$ equipped with a lift of the circle action.  The quotient of the holomorphic bundle upstairs by $U(1)$ produces the pair of holomorphic bundles over $\Sigma$ downstairs.  

Choose the 4-manifold $X$ so that the fibre bundle $X-\{\hat{p}_1,...,\hat{p}_k\}/U(1)\cong\Sigma\times I-\{y_1,...,y_k\}$ restricts to the trivial $U(1)$ bundle above $\Sigma\times\{0\}$, and to the $U(1)$ bundle $\mathcal{O}(z_1,...,z_k)$ above $\Sigma\times\{1\}$.  Choose $V$ harmonic on $\Sigma\times I-\{y_1,...,y_k\}$ with $V\sim 1/2r$ near each $y_j$ ($r$ is the distance from $y_j$) and 
\begin{equation}  \label{eq:harmboun}
\partial V/\partial t|_{\Sigma\times\{0\}}=0,\quad V|_{\Sigma\times\{1\}}=1.
\end{equation}  
Since $\left[\frac{1}{2\pi}*dV\right]$ is an integral cohomology class, there exists a connection $d\theta+a$ on the circle bundle 
\[ X-\{\hat{y}_1,...,\hat{y}_k\}\to\Sigma\times I-\{y_1,...,y_k\}\] 
with curvature $da=*dV$.  The connection corresponds to a $U(1)$-invariant 1-form $\omega$ on $X$.  The Gibbons-Hawking metric 
\[ g=Vds^2+V^{-1}\omega^2\]
extends over each $\hat{p}_j$ to give a smooth $U(1)$-invariant metric on $X$.  (We have abused terminology and identified $V$ and $ds^2$ with their pull-backs upstairs.)  The harmonic function $V$ uniquely determines a $U(1)$ monopole $(A_V,\Phi_V)$ with $A_V$ trivial on $\Sigma\times\{0\}$ and $\Phi_V=i(V-1)$.  The restriction of $A_V$ to $\Sigma\times\{0\}$, respectively $\Sigma\times\{1\}$,  determines two holomorphic line bundles over $\Sigma$, $L_0=\mathcal{O}$, the trivial line bundle, and $L_1=\mathcal{O}(z_1,...,z_k)$.

There is a complex structure on $X-\{\hat{p}_1,...,\hat{p}_k\}$ given by 
\[ J\cdot d\bar{z}=id\bar{z},\quad J\cdot Vdt=\omega.\]  
Here we identify $d\bar{z}$, $dt$ and $V$ with their pull-backs upstairs on $X$. 
\begin{lemma}
The complex structure $J$ is integrable and extends to $X$ .  With respect to $J$, the $U(1)$-invariant metric g is K\"ahler.
\end{lemma}
\begin{proof}
We must show that $J$ is integrable on $X-\{\hat{p}_1,...,\hat{p}_k\}$, and that the K\"ahler form of $g$ is closed.  It will follow that $J$ is horizontal, so it uniquely extends to $X$ via parallel transport.  The proof follows LeBrun \cite{LeBExp}.

The Newlander-Nirenberg condition for the integrability of $J$ is 
\[ [T_X^{1,0},T_X^{1,0}]\subset T_X^{1,0}\] 
where $T_X^{1,0}$ is the $+i$ eigenspace of $J$ on $TX\otimes\bc$, equivalently defined as the kernel of $(d\bar{z},Vdt-i\omega)$.

Let $u$ and $v$ be vector fields taking values in $T_X^{1,0}$.  Then $\eta(u)=0=\eta(v)$ for any 1-form $\eta$ in the span of $d\bar{z}$ and $Vdt-i\omega$.  The integrability condition for $J$ is equivalent to the vanishing $\eta[u,v]=0$.  This follows from the identity \[\eta[u,v]=u(\eta(v))-v(\eta(u))-d\eta(u,v)\]  together with knowledge of the exterior derivatives of $d\bar{z}$ and $Vdt-i\omega$.  We have $d(d\bar{z})=0$ and:
\begin{eqnarray*} 
d(Vdt-i\omega)&=&dV\wedge dt-id\omega\\
&=&(V_zdz+V_{\bar{z}}d\bar{z})\wedge dt-i(-iV_zdz\wedge dt+iV_{\bar{z}}d\bar{z}\wedge dt+2i\rho V_tdz\wedge d\bar{z})\\ 
&=&2V_{\bar{z}}d\bar{z}\wedge dt +2\rho V_tdz\wedge d\bar{z}\\
&=&d\bar{z}\wedge(2V_{\bar{z}}dt-2\rho V_t dz)
\end{eqnarray*}
which vanishes on $T_X^{1,0}$ as required.  (Although $\omega$ is not the pull-back of a 1-form downstairs, $d\omega$ is the pull-back of the 2-form $*dV$.)

The K\"{a}hler form of $g$ is
\[\Omega=2i\rho Vdz\wedge d\bar{z}+dt\wedge\omega\] 
and the $t$ independence of $\rho$ guarantees that it is closed
\[d\Omega=2i\rho dV\wedge dz\wedge d\bar{z}-dt\wedge d\omega=2i\rho V_tdt\wedge dz\wedge d\bar{z}-2i\rho V_tdt\wedge dz\wedge d\bar{z}=0.\]  
\end{proof}
The $U(1)$ action on $X$ extends to a local $\bc^*$ action---i.e. an action of a neighbourhood of $\{|z|=1\}$ in $\bc^*$---on $X$ (and the local action on the boundary of $X$ has the further restriction that only $\{|z|\geq 1\}$ acts.)  The $\bc^*$ action in a neighbourhood of a fixed point $\hat{y_j}$ is modeled by the $\bc^*$ action on $\bc^2$ near $(0,0)$
\begin{equation}   \label{eq:local}
\lambda\cdot(x,y)=(\lambda x,\lambda^{-1}y).
\end{equation}
Hence, as a complex manifold, $X$ can be realised as a submanifold of the algebraic variety
\[ \tilde{X}=\{(x,y)\in L_0\oplus L_1|\ xy=s\}\] 
where $s\in H^0(L_1)$ is a section that vanishes at $z_1,...,z_k$, and $L_1=\mathcal{O}(z_1,...,z_k)$ and $L_0=\mathcal{O}$ are the holomorphic line bundles over $\Sigma$ determined by the $U(1)$ connection $A_V$.    The algebraic surface $\tilde{X}$ maps to $\Sigma$ with generic fibres $\bc^*$ and nodal singularities above each $z_j$.  Kapustin and Witten showed that the complex structure on $X$ is independent of the $t$ coordinate of each $y_j=(z_j,t_j)\in\Sigma\times I$ so it is not surprising that this construction depends only on the points $z_j$.

The trivial $G$-bundle $E$ over $\Sigma\times I$ pulls back to a $U(1)$-invariant bundle $\tilde{E}$ over $X$ equipped with a lift of the $U(1)$-action.  A singular monopole $(A,\Phi)$ on $E$ lifts to a $U(1)$-invariant anti-self-dual connection $\tilde{A}$ on $\tilde{E}$, i.e. $F_{\tilde{A}}=-*F_{\tilde{A}}$, which defines a holomorphic structure on $\tilde{E}$.  The local $\bc^*$ action on $X$ lifts to a local $\bc^*$ action on $\tilde{E}$.  The bundle $\tilde{E}$ is equipped with a framing over any basepoint above $\Sigma\times\{0\}$. The framing upstairs is the lift of a framing of $E$ over the image downstairs of the basepoint, that extends to a framing over all of $\Sigma\times\{0\}$ using the flat connection $A|_{\Sigma\times\{0\}}$ and hence is equivalent to requiring gauge transformations to be the identity over $\Sigma\times\{0\}$.

The triviality of the $U(1)$ connection $A_V$ on $\Sigma\times\{0\}$ means that any constant section $\Sigma\times\{0\}$ of the trivial $U(1)$ bundle is a holomorphic section $\Sigma\to X$.    We denote the image holomorphic curve of this section by $\Sigma_0$.  (There is a circle of such holomorphic sections.) 
The restriction of $\tilde{E}$ to $\Sigma_0$ has a canonical holomorphic trivialisation.  This is because any other trivialisation differs by a holomorphic map $\Sigma_0\to G^c$ that is the identity over the basepoint, and hence the identity everywhere by the compactness of $\Sigma_0$.  This extends to a canonical $\bc^*$-invariant trivialisation of $\tilde{E}$ over the dense subset of $X$ given by the orbit of the local $\bc^*$-action on $X$.

Given a disk $D\subset\bc$ containing 0 consider two embeddings of $D$ into a neighbourhood of a fixed point $\hat{y}_j$ with images $D_+=\{(\epsilon,y)\in\bc^2|y\in D\}$  and $D_-=\{(x,\epsilon)\in\bc^2|x\in D\}$ for small $\epsilon$ with respect to the local coordinates around $\hat{y}_j$ given in (\ref{eq:local}).  The $\bc^*$ action defines an isomorphism between $D_+$, $D_-$ and a disk $D_{z_j}\subset\Sigma_0$ containing $z_j$.  It lifts to define an isomorphism
\[h:\tilde{E}|_{D_+-\{0\}}\to\tilde{E}|_{D_--\{0\}}.\]  
Moreover, the canonical $\bc^*$-invariant trivialisation of $\tilde{E}$ over a dense subset of $X$ determines a trivialisation over $D_+$ so $h$ gives a map $D-\{0\}\to G^c$ well-defined up to changes of trivialisation of $\tilde{E}$ over $D_-$ thus determining an element $\gamma_j\in LG^c/L^+G^c$.  

Conversely the construction above shows that elements $\gamma_j\in LG^c/L^+G^c$ and disks $D_{z_j}\subset\Sigma_0$, for $j=1,...,k$, determine a $\bc^*$-invariant holomorphic bundle $\tilde{E}$ over $X$.  Thus the question of existence and uniqueness of singular monopoles on the trivial bundle over $\Sigma\times I$ has been translated to a question of existence and uniqueness of $U(1)$ invariant instantons on a $\bc^*$-invariant holomorphic bundle $\tilde{E}$ over the K\"ahler 4-manifold $X$.   This problem was solved by Donaldson who proved in \cite{DonAnt} that a stable holomorphic bundle over a closed K\"ahler 4-manifold $X$ possesses a unique instanton.  Furthermore, Donaldson showed in \cite{DonBou} that if $X$ has boundary then no stability condition is required on the holomorphic bundle.

Donaldson's proof uses the fact that a Hermitian metric $\tilde{H}$ on a holomorphic bundle $\tilde{E}$ over a complex surface determines a connection $\tilde{A}=\tilde{H}^{-1}\partial\tilde{H}$ on the bundle.  The Hermitian metric $H$ used in Section~\ref{sec:uniq} lifts to a $U(1)$ invariant Hermitian metric $\tilde{H}$ on $\tilde{E}$ over $X$.  Prescribe $\tilde{H}|_{\partial X_0}=I$ and $\partial\tilde{H}/\partial n|_{\partial X_1}=0$ where $\partial X=\partial X_0\cup\partial X_1$ for $\partial X_i$ corresponding to the boundary $\Sigma\times\{i\}$, $i=0,1$ downstairs.  These correspond to the boundary conditions (\ref{eq:bou}) for the singular monopole downstairs.
\begin{thm}[Donaldson \cite{DonBou}]
Let $\tilde{E}$ be a holomorphic vector bundle over the compact K\"ahler manifold $X$, with non-empty boundary $\partial X$.  For any Hermitian metric $f$ on the restriction of $\tilde{E}$ to $\partial X$ there is a unique Hermitian metric $\tilde{H}$ on $\tilde{E}$ such that

(i) $\tilde{H}=f$ over $\partial X$,

(ii) $F_{\tilde{H}}$ satisfies the Hermitian Yang-Mills equation,\\
where $F_{\tilde{H}}$ is the curvature of the connection associated to $\tilde{H}$.
\end{thm}
This theorem gives a solution to a Dirichlet boundary problem however we need an adjusted version of this that allows a mixture of Dirichlet and Neumann boundary conditions.  Such a Neumann version is implicit in Simpson's work \cite{SimCon} --- there exists a unique Hermitian Yang-Mills metric $\tilde{H}$ (equivalently $\tilde{A}$ is anti-self-dual with respect to the  K\"ahler metrc) for any specified boundary value $\tilde{H}|_{\partial X}=H_0$.   Uniqueness of the solution together with $U(1)$-invariance of the boundary conditions guarantees that the solution is $U(1)$ invariant and hence descends to the required singular monopole.   Thus we get a second proof of existence and uniqueness of singular monopoles presented in Corollary~\ref{th:surj} and Proposition~\ref{th:1to1}.  The direct proofs of Corollary~\ref{th:surj} and Proposition~\ref{th:1to1} are essentially a dimensional reduction of the proof given by Donaldson in \cite{DonBou} that avoid the choice of a (singular) harmonic function $V:\Sigma\times I\to\br$ which determines the K\"ahler metric.


\begin{thebibliography}{99}

\bibitem{AHiGeo} Atiyah, Michael and Hitchin, Nigel 
\emph{The geometry and dynamics of magnetic monopoles.} 
Princeton University Press, Princeton, NJ, 1988.

\bibitem{DonAnt} Donaldson, Simon
\emph{Anti-self-dual Yang-Mills connections over complex algebraic surfaces and stable vector bundles.} 
Proc. London Math. Soc. {\bf 50} (1985) 1-26. 

\bibitem{DonBou} Donaldson, Simon
\emph{Boundary value problems for Yang-Mills fields.} 
J. Geom. Phys. {\bf 8} (1992), 89-122. 

\bibitem{DonNah} Donaldson, Simon
\emph{Nahm's equations and the classification of monopoles.} 
Comm. Math. Phys. {\bf 96} (1984), 387--407.

\bibitem{FBeVer}   Frenkel, Edward and Ben-Zvi, David 
\emph{Vertex algebras and algebraic curves. Second edition.}
 Mathematical Surveys and Monographs {\bf  88}. AMS, Providence, RI, 2004.
 
\bibitem{GraHol} Grauert, Hans 
\emph{Holomorphe Funktionen mit Werten in komplexen Lieschen Gruppen.}  Math. Ann. {\bf 133} (1957), 450-472.

\bibitem{JarCon} Jarvis, Stuart
\emph{Construction of Euclidean monopoles.}
Proc. London Math. Soc. (3) {\bf 77} (1998), 193--214.

\bibitem{JNoCom} Jarvis, Stuart and Norbury, Paul
\emph{Compactification of hyperbolic monopoles.} 
Nonlinearity {\bf 10} (1997) 1073-1092. 

\bibitem{KWiEle}  Kapustin, Anton  and Witten, Edward 
\emph{Electric-Magnetic Duality And The Geometric Langlands Program.} 
Commun. Number Theory Phys. {\bf 1} (2007), 1-236. 
 
\bibitem{KroMon} Kronheimer, Peter
\emph{Monopoles and Taub-NUT metrics.}
MSc thesis, Oxford (1985). 

\bibitem{LeBExp} LeBrun, Claude
\emph{Explicit dual metrics on $\bc\bp_2\#...\#\bc\bp_2$.}
J. Diff. Geom. {\bf 34} (1991), 223-253.

\bibitem{PauMon} Pauly, Marc 
\emph{Monopole moduli spaces for compact 3-manifolds.} 
Math. Ann. {\bf 311} (1998), 125-146.

\bibitem{PSeLoo} Pressley, Andrew and Segal, Graeme 
\emph{Loop groups.} 
Oxford Mathematical Monographs. Oxford University Press, 1986.

\bibitem{SimCon}  Simpson, C.T.  
\emph{Constructing variations of Hodge structure using Yang-Mills 
theory, and applications to uniformisation.} 
J. Am. Math. Soc. {\bf 1} (1989) 867-918. 

\bibitem{TayPar} Taylor, Michael E.
\emph{Partial differential equations I.} 
Springer-Verlag, New York, 1996.


\end{thebibliography}
\end{document}